\documentclass[final,3p]{elsarticle}

\usepackage{amssymb}
\usepackage{amsmath,amsfonts}
\usepackage{mathtools}
\usepackage{xcolor}
\usepackage{booktabs, multirow} 
\usepackage{soul}
\usepackage{graphicx}
\usepackage{caption}
\usepackage{subcaption}
\usepackage{float}
\usepackage{gensymb}
\usepackage{ulem}

\usepackage{cleveref}

\newcommand{\ssigma}{\boldsymbol{\sigma}}

\newcommand{\eref}[1]{Equation~(\ref{#1})}

\newcommand{\fref}[1]{Figure~\ref{#1}}
\newcommand{\frefs}[1]{Figures~\ref{#1}}
\newcommand{\tref}[1]{Table~\ref{#1}}
\newcommand{\sref}[1]{Section~\ref{#1}}

\biboptions{sort&compress}

\begin{document}

\begin{frontmatter}


%
\title{On the fractional transversely isotropic  functionally graded nature of soft biological tissues}
%

\author[IITM]{G. Sachin}
\ead{me21d011@smail.iitm.ac.in}
\author[IITM]{S. Natarajan}
\ead{snatarajan@smail.iitm.ac.in}
\author[OBU,UO]{O. Barrera\corref{cor1}}
\ead{olga.barrera@eng.ox.ac.uk, +44 7585041147.}

\cortext[cor1]{Corresponding author}

\address[IITM]{Department of Mechanical Engineering, Indian Institute of Technology Madras, Chennai, 600036, India.}
\address[OBU]{School of Engineering, Computing and Mathematics, Oxford Brookes University, Headington, Oxford OX3 0BP, United Kingdom.}
\address[UO]{Department of Engineering Science, University of Oxford, Parks Road, OX1 3PJ, Oxford, United Kingdom.}

\begin{abstract}
This paper focuses on the origin of the poroelastic anisotropic behaviour of the meniscal tissue and its spatially varying properties. We present confined compression creep test results on samples extracted from three parts of the tissue (Central body, Anterior horn and Posterior horn) in three orientations (Circumferential, Radial and Vertical). We show that a poroelastic model in which the fluid flow evolution is ruled by non-integer order operators (fractional Darcy's law) provides accurate agreement with the experimental creep data. The model is validated against two additional sets of experimental data: stress relaxation and fluid loss during the consolidation process measured as weight reduction. Results show that the meniscus can be considered as a transversely isotropic poroelastic material. This behaviour is  due to the fluid flow rate being about three times higher in the circumferential direction than in the radial and vertical directions in the body region of the meniscus.
In the anterior horn, the elastic properties are transversely isotropic, with the aggregate modulus higher in the radial direction than in the circumferential and vertical directions.
The 3D fractional poroelastic model is implemented in finite element software and quantities  such as flux of interstitial fluid during the consolidation process, a non-trivial experimental measure, are determined.

\end{abstract}

\end{frontmatter}

\section{Introduction}
\begin{figure}[H]
            \centering
            \begin{subfigure}{.4\textwidth}
                \centering
                  \includegraphics[width=\linewidth]{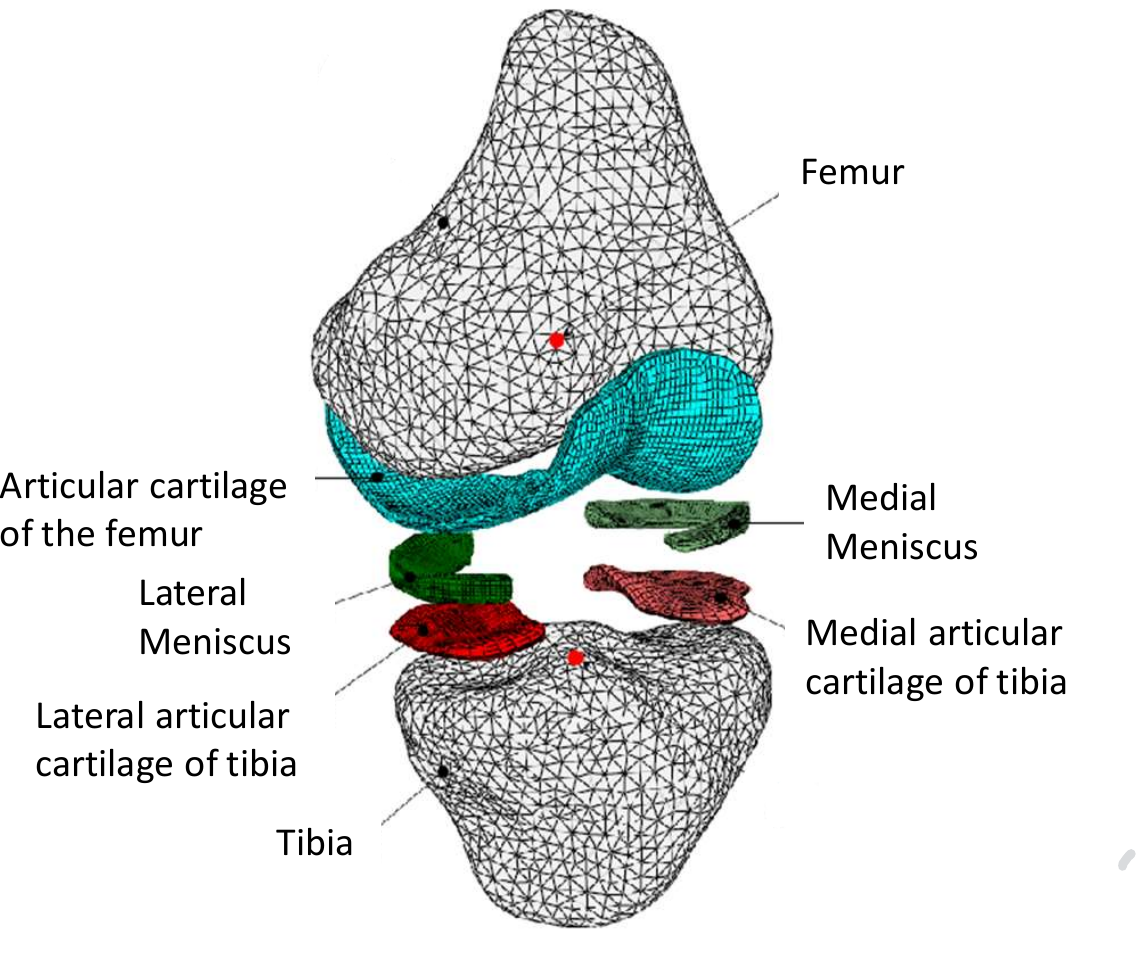}
                \caption{}
                \label{fig:KneeJoint}
            \end{subfigure}%
        \begin{subfigure}{.6\textwidth}
            \centering
            \includegraphics[width=\linewidth]{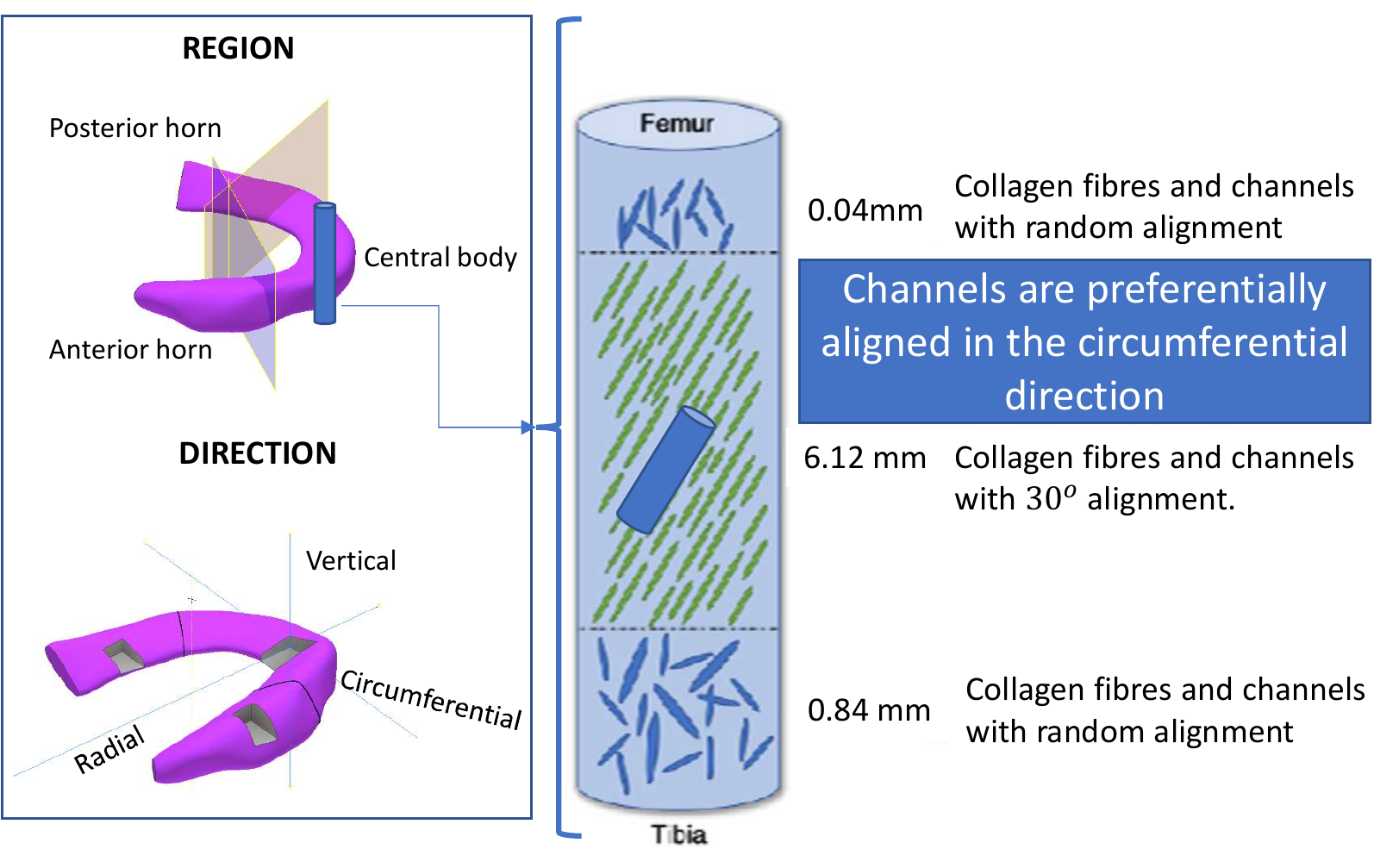}
            \caption{}
            \label{fig:RegionsnDirections}
        \end{subfigure}\\%
            \begin{subfigure}{.4\textwidth}
            \centering
            \includegraphics[width=\linewidth]{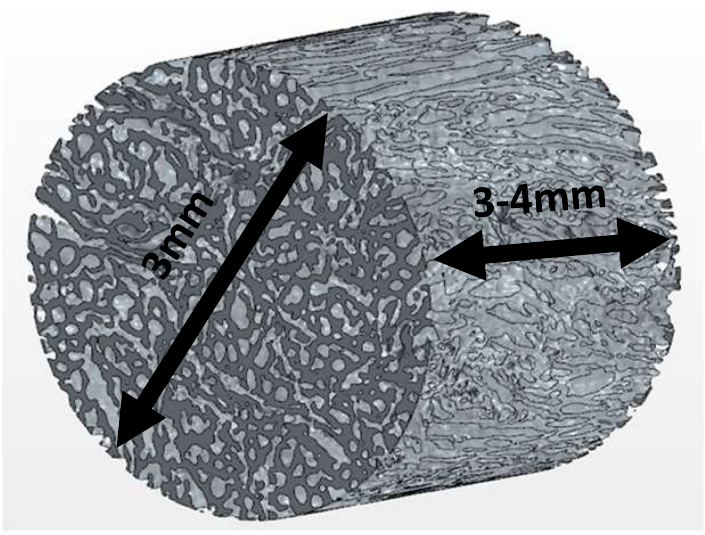}
            \caption{}
            \label{fig:porousCylinder}
        \end{subfigure}%
        \begin{subfigure}{.6\textwidth}
            \centering
            \includegraphics[width=\linewidth]{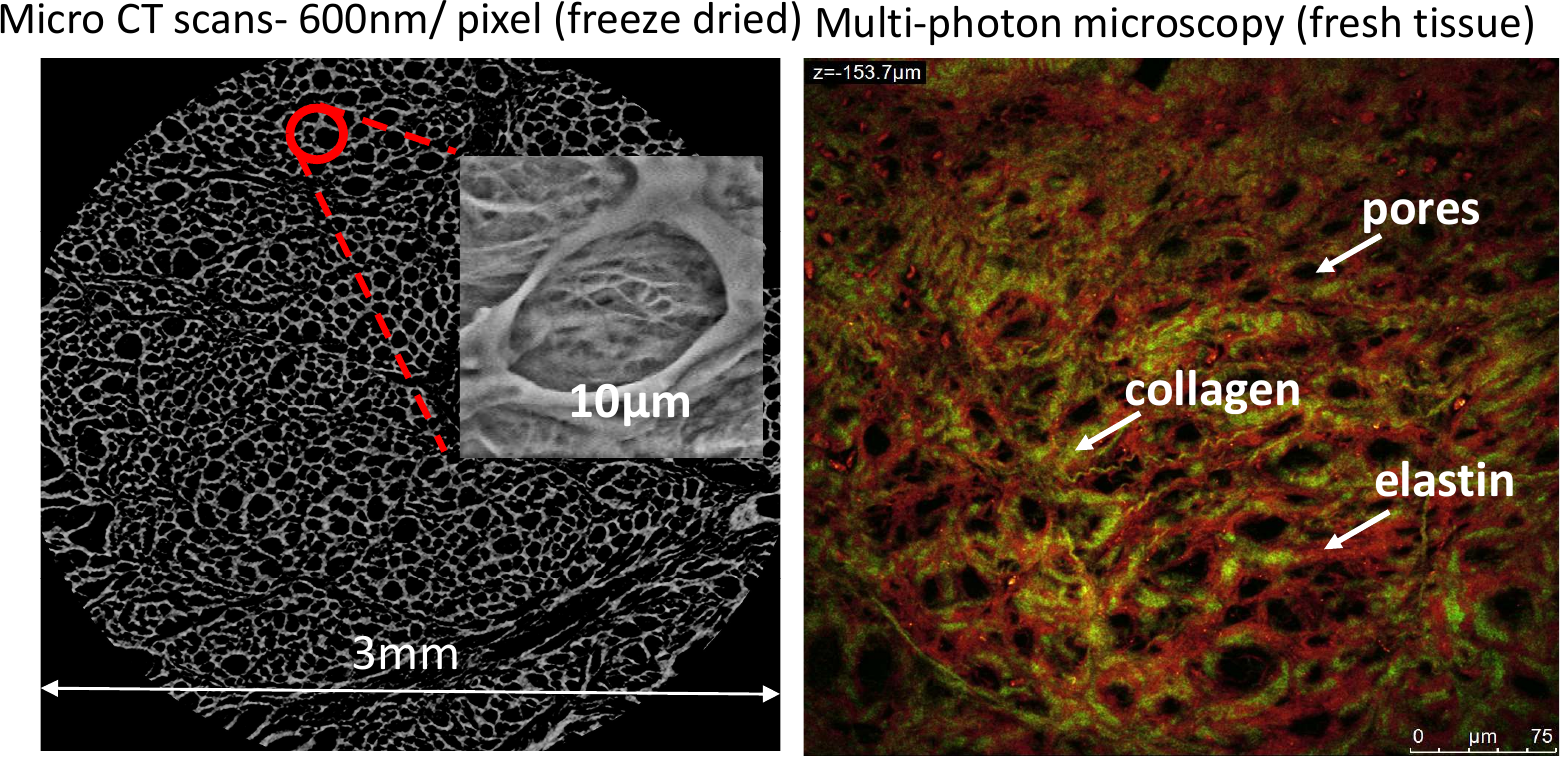}
            \caption{}
            \label{fig:CTscans}
        \end{subfigure}\\%
        \caption{(a) Lateral and Medial menisci in the knee joint, (b) Regions and the directions in the meniscus, (c) 3D porous microstructure of meniscus cylindrical sample (d)Micro CT scan and Multi-photon microscopy of meniscus tissue.}
        \label{fig:knee}
    \end{figure}

Biological soft tissues, such as the meniscal tissue (\fref{fig:knee}), exhibit a hierarchical porous solid matrix and an interstitial fluid flowing into the pores \citep{Agustoni1, Maritz1, Bonomo2019, VetrietAl, Bulle2021}. The overall mechanical behaviour 
depends not only on the solid matrix deformation but also on the movement of the fluid in and out of the pores during the deformation. Furthermore, the internal architecture of the structure constituting the solid matrix and collagen determines the anisotropic behaviour of both the elastic and the transport properties \citep{Ateshian2010}. The poroelastic theory is employed to model the response of biphasic tissue and better understand the anisotropic nature of the meniscus on the functioning of the knee.


The menisci situated between the femoral and tibial cartilage of the knee joint (\fref{fig:KneeJoint}) have three main regions, viz., anterior horn, central body and posterior horn (\fref{fig:RegionsnDirections}). They perform several functions like load bearing, lubrication, energy absorption and stability \cite{Kurosawa, Shrive}. In light of our recent experimental, theoretical and numerical studies \citep{Agustoni1, Maritz1, Biaxial, kneeFEM, imagedriven}, we understand that each region has different functions. 
In \citep{Maritz1, Agustoni1}, the unique structure of the body region that fulfils load-bearing and energy absorption capabilities is discussed.
This is achieved through a sandwich-like structure with two thin, stiff outside layers which take the load-bearing role and a thick, softer internal layer which acts as an energy absorption element, an effective natural damper \citep{Maritz1,imagedriven}. The architectural arrangement of the solid collagen-based matrix appears to be different and gradually changing in the outside and the internal layers. The outside layers show a dense distribution of collagen fibres with random orientation and low mean diameter. The elastic modulus is one order of magnitude higher than the internal layer. Whereas the internal layer of the meniscus relies on a hierarchical anisotropic network of collagen channels mainly aligned along the circumferential direction in which fluid flows during the physiological deformation process (\fref{fig:knee}). In \citep{Bulle2021}, we have shown that the permeability changes from the anterior/posterior horns to the body region, with the body region being more permeable. Here, we show that the permeability tensor in this region is also transversely isotropic and that the permeability in the circumferential direction (the preferential direction of the collagen channel) is higher than the other two directions (\fref{fig:knee}). 
Fluid flowing inside these channels determines the time-dependent behaviour of the tissue. As the tissue deforms under the action of loads, these channels change morphology, which results in change in permeability as a function of channels dimension, porosity, and tortuosity. Thus, making the permeability tensor not only depending on the pressure but also varying spatially and temporally because of the local variation of applied pressure gradients.
Recent finite element studies of the human knee have shown that anterior and posterior horns might have a stability function \citep{ELMUKASHFI2022,kneeFEM} as they deform to allow the meniscus to be extruded 
so that the body region has a higher contact area with the cartilage and can then carry most of the load. 

Appropriate modelling of the fluid flow behaviour in the different portions of the meniscal tissue when subject to a range of loading conditions is essential to gain insight into the biomechanical function of the tissue and design appropriate artificial tissue substitutes. Current barriers to the clinical and functional outcomes
of the to-date devices \citep{vanKampen2013} are that the transport and structure-properties relationship of the native tissues is
still not well understood. In particular, the relationship between the ``micro/nano''-scopic composition and structure of the meniscus and its coarse-scale anisotropic behaviour structure is still an open question \citep{FoxAJ}.

Experimental tests, such as confined compression tests (relaxation and/or creep) are currently used to characterize the material parameters such as elastic modulus at equilibrium and permeability. The main models used to identify these two parameters are based on biphasic theory.
In \citep{Soltz1998}, the authors examined, through a confined compression test, the effect of fluid in cartilage is bearing up to 90\% of the applied load. They used a linear biphasic model to find out the material properties in terms of the aggregate modulus and permeability. In \citep{Sweigart2003}, authors performed a creep indentation test on different regions of the meniscus of different animal species. The different material properties such as the aggregate modulus, Poisson's ratio, shear modulus and permeability were estimated by using the biphasic model to fit the experimental data.


Our recent experimental findings on human meniscal samples show that the biphasic model was not sufficient to reproduce the observed experimental behaviour \citep{Bulle2021}. To address this, \cite{Barrera2021} proposed a poroelastic model wherein the pore pressure diffusion equation is derived by adopting a modified version of Darcy's law involving fractional derivatives. It follows that the hydraulic permeability becomes ``anomalous" and the rate of fluid flow is governed by the order of the derivative.
It is thus clear that the anomalous hydraulic permeability, the order of derivative tensors and their anisotropic nature, which govern the transport of fluid within the complex porous structure of the meniscus, are key properties of the meniscal tissue. Though a finite element implementation of the fractional poroelastic model is reported in \citep{Barrera2021}, however, the model is restricted to estimating only the pore pressure field.


The purpose of this paper is to present for the first time creep data of confined compression tests performed in different regions of the meniscal tissue (anterior-horn, central-body and posterior-horn) and in different directions (radial, circumferential and vertical). We adopt the fractional poroelastic model presented in \citep{Barrera2021, Bulle2021} to fit the experimental data to recover the aggregate modulus, anomalous permeability and order of derivative and their dependence on the directions (Aggregate modulus, anomalous permeability and order of derivative tensors). For the first time, we validate the model in terms of Stress, weight loss and displacement fields by using the parameters recovered from fitting the creep curve to match two additional sets of data for the same samples, i.e., relaxation and weight loss data. The fractional poroelastic model is implemented in commercial software Abaqus using UMATHT and UMAT. The numerical simulation results show a good agreement with analytical solutions regarding both displacements and pore pressure fields which were lacking in the literature. 
 This enabled us to run consolidation tests numerically and calculate fields such as the flux of interstitial fluid during the consolidation process, which is a difficult experimental measure to obtain. Additionally, we perform statistical analysis of the parameters using ANOVA 
method to compare the variation of the material parameters in different directions. 

\section{Material and Methods}
\subsection{Confined consolidation creep tests}
\label{sec:ConfinedCompression}
Menisci were harvested from patients undergoing total knee arthroplasty under ethically approved protocol as reported in \citep{Bulle2021}. Samples labelled as ``degraded" by the gross investigation of the surgeon were discarded. Menisci were thawed in Phosphate-buffered solution (PBS) for thirty minutes to recover their hydration. A total of 29 meniscal test samples were harvested. The samples were cylindrical in shape with a 3mm diameter and 3-4 mm in height. Samples were extracted from the central body, anterior horn and posterior horn regions along each principal direction, i.e. circumferential, radial and vertical (\fref{fig:knee}).
\begin{figure}[H]
\centering
\begin{subfigure}{.5\textwidth}
\centering
\includegraphics{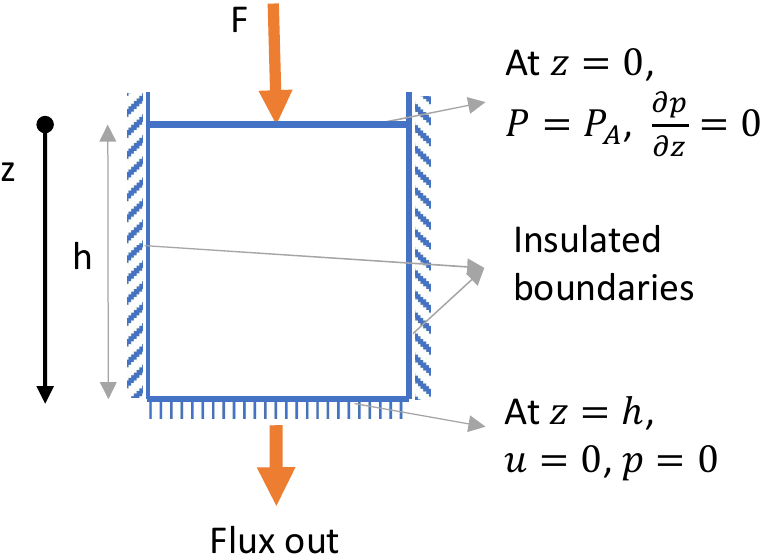}
\caption{}
\label{sf:Schematic_Exp}
\end{subfigure}\\%
\begin{subfigure}{.5\textwidth}
\centering
\includegraphics[width=0.9\linewidth]{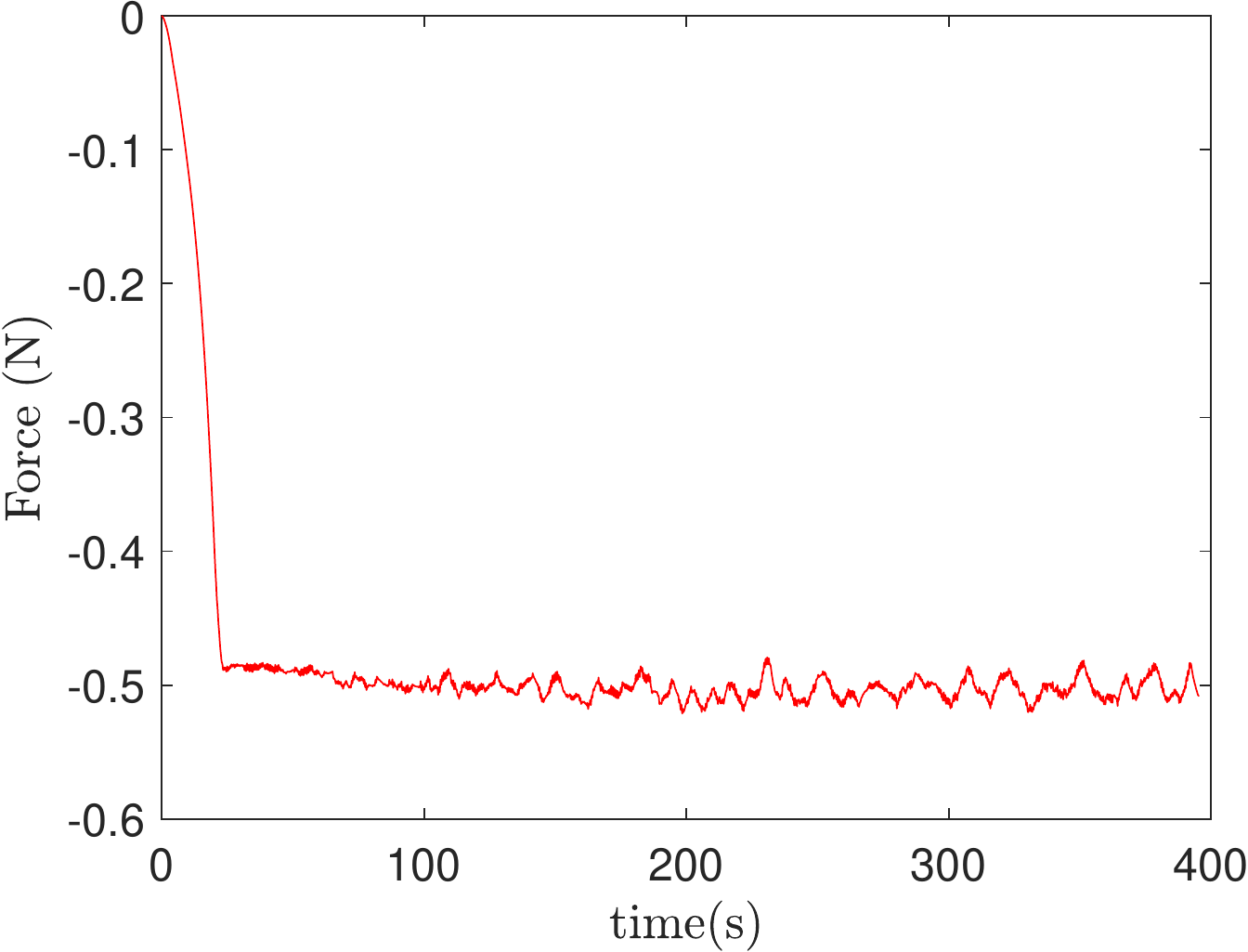}
\caption{}
\label{sf:ExpForce}
\end{subfigure}%
\begin{subfigure}{.5\textwidth}
\centering
\includegraphics[width=0.9\linewidth]{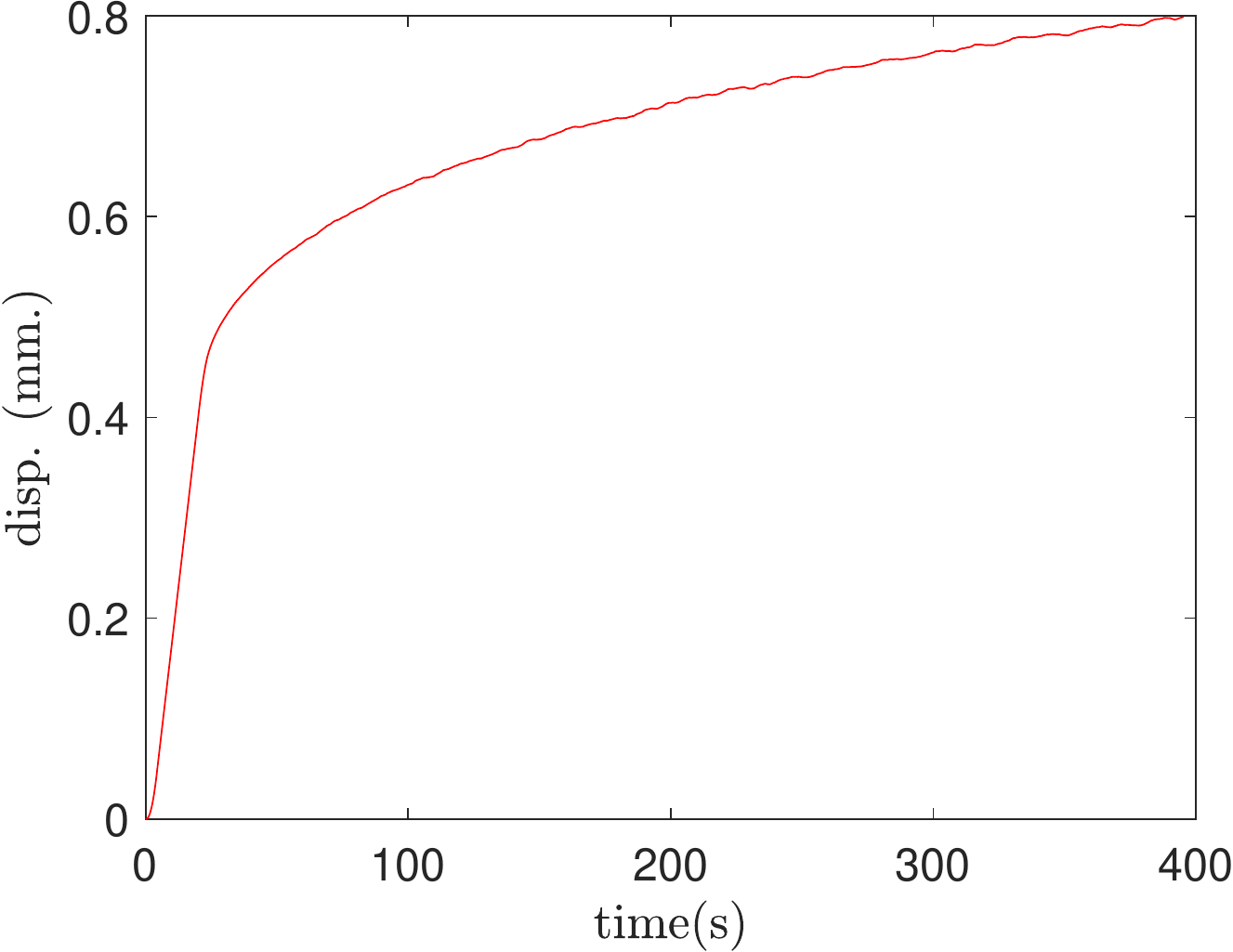}
\caption{}
\label{sf:ExpDisp}
\end{subfigure}%
\caption{(a) Schematic representation of confined compression test, (b) applied load for confined compression creep test and (c) measured displacement.}
\label{fig:boundcond}
\end{figure}


The confined compression test set-up consisted of a chamber with transparent and insulated walls and a porous plate at the bottom(\fref{sf:Schematic_Exp}). The sample to be tested is confined inside the chamber. A load of 0.5 N (corresponding to a stress of about 0.07 MPa which is a physiological value for menisci \cite{SEITZ201368}) at the top of the chamber is applied. It is applied through a piston with a stage velocity of 0.3 \% height of the sample till the load 0f 0.5 N is reached and then the load is kept constant for 400s as shown in \fref{fig:boundcond}. The applied force and the obtained displacement that is measured at the top surface are shown in \frefs{sf:ExpForce}-\ref{sf:ExpDisp}, respectively. 

\subsection{Fractional consolidation - governing equations \& finite element implementation}
\paragraph{Governing equations}
\label{sec:governing eqns}
Consider a poroelastic material occupying $\Omega \subset \mathbb{R}^d$, bounded by $\partial \Omega \subset \mathbb{R}^{d-1}$, where $d=1,2,3$. The boundary accommodates the following decomposition: $\partial \Omega = \partial \Omega_N \cup \partial \Omega_D$ and $\partial \Omega_N \cap \partial \Omega_D = \emptyset$, where Neumann and Dirichlet boundary conditions are applied on $\partial \Omega_N$ and $\partial \Omega_D$, respectively. In the absence of body force and inertia, the governing differential for a poroelastic material is given by:
\begin{subequations}
\begin{align}
    \boldsymbol{\nabla} \cdot \ssigma &= \mathbf{0} \\
    \dfrac{\partial \zeta}{\partial t} + \boldsymbol{\nabla} \cdot \mathbf{J}_{p} &= 0
\end{align}
\label{eqn:strongform}
\end{subequations}
supplemented with the following boundary conditions: $\mathbf{u}=\hat{\mathbf{u}},~{\rm and}~p = \hat{p}$ on $\partial \Omega_D$ and $\ssigma \cdot \mathbf{n} = F,~{\rm and}~\mathbf{J}_p \cdot \mathbf{n}= 0$ on $\partial \Omega_N$. In \eref{eqn:strongform}, $\zeta$ is the variation of the fluid content, $p$ is the pore pressure, $\mathbf{J}_{p}$ is the fluid flux and $\ssigma$ is the stress. The relation between the stresses, displacements and pressure; and between the flux and the pressure is given by:
\begin{subequations}
\begin{align}
    \ssigma &= 2G\boldsymbol{\epsilon} + \lambda \mathrm{Tr}(\boldsymbol{\epsilon})\textbf{I} - \alpha p \textbf{I} \\
    \mathbf{J}_p &= -\lambda_{\beta} \prescript{}{0}D_{t}^{\beta} \boldsymbol{\nabla} p \label{eq:FluidFlux}
\end{align}
\label{eqn:governeqn}
\end{subequations}
where $\lambda = \dfrac{(3K-2G)}{3}$ is Lam\'e's constant, $G, K$ are the shear and the bulk modulus, respectively, $\alpha$ is the Biot coefficient and $\prescript{}{0} D_{t}^{\beta}$ is the fractional derivative with respect to time of order $\beta$. It is important to note that in case of $\beta=$ 0, $\lambda_{\beta}$ is with dimension of $[L]^{4}[F]^{-1}[T]^{-1}$. In the case of $\beta \neq 0$, $\lambda_{\beta}$ has dimension of $[L]^{4}[F]^{-1}[T]^{\beta-1}$. The pore pressure diffusion equation (c.f. \eref{eq:FluidFlux}), written in terms of dilatational strain, $\epsilon_{d}$, is given by:
\begin{equation}
\dfrac{\partial p}{\partial t} = \dfrac{K_{u}-K}{\alpha^{2}} \lambda_{\beta} \prescript{}{0}D_{t}^{\beta} \boldsymbol{\nabla}^{2}p - \dfrac{K_{u} - K}{\alpha}\dfrac{\partial \epsilon_{d}}{\partial t}
\label{eq:ppDiffusion}
\end{equation} 
where $K_u$ is the undrained bulk modulus.

During confined compression tests, a cylinder of meniscal tissue is inserted inside a cylindrical chamber with a porous base, see \fref{sf:Schematic_Exp} for geometry and boundary conditions. Upon applying pressure at the top of the sample, the fluid contained inside the sample will flow out of the base. This can be modelled as a one-dimensional version of \eref{eqn:governeqn} as the only non-zero components of strain will be $\epsilon_{zz}$. \eref{eqn:governeqn} will be modified to~\citep{Barrera2021}: 

\begin{subequations}
\begin{align}
\left(K+\frac{4G}{3}\right) \dfrac{\partial \epsilon_{zz}}{\partial z} - \alpha \dfrac{\partial p}{\partial z}=0  
\label{eq:Strnpp} \\
\dfrac{\partial p}{\partial t} = \dfrac{K_u - K}{\alpha^2}\lambda_\beta  \prescript{}{0}D_t^\beta \dfrac{\partial^2p}{\partial z^2} - \dfrac{K_u - K}{\alpha}\dfrac{\partial \epsilon}{\partial t}
\end{align}
\end{subequations}
and the corresponding boundary conditions are, $\forall~ 0 \leq z \leq h$, 
$ p(z,t=0)=\gamma P_{A};~P(z=0,t)=P_{A},~p(z=h,t)=0,~
\dfrac{\partial p}{\partial z}\Big|_{z=0}=0$ and $ 
u(z=h,t)=0.$
where $h$ is the initial height of the sample. A constant compressive stress, $\sigma_{zz}(z=0,t)=- P_{A}$ in the $z-$ direction is applied to the cylinder at $z=0$.
The initial pore pressure is derived assuming that the material is under undrained conditions at the initial time of loading, given by:
\begin{equation}
p(z,t=0)= P_{A}\frac{3(K_{u}-K)}{\alpha\left(4G+3K_{u}\right)}
\label{eq:9}
\end{equation}
It is shown that the analytical solution of the pore pressure is given by~\citep{Barrera2021}:
\begin{equation}
p(z,t)= P_{A}\gamma \sum \limits_{n=1,3}^\infty  E_{1-\beta,1} \left(-{\frac{n^{2}\pi^{2} \bar{\lambda} t^{1-\beta}}{4h^{2}}}\right) c_{n}\cos \frac{n \pi z}{2h}
\label{eq:PorepSoln}
\end{equation}
where $E_{1-\beta,1}$ is the Mittag-Leffler function and
\begin{equation}
\begin{aligned}
\gamma &= \dfrac{3\left(K_{u}-K\right)}{\alpha\left(4G+3K_{u}\right)}; \qquad
c_{n} = \dfrac{4}{n \pi} \big(-1\big)^{\frac{n-1}{2}}\\
\bar{\lambda} &= \lambda_{\beta}\dfrac{\left(4G+3K\right)\left(K_{u}-K\right)}{\alpha^{2}\left(4G+3K_{u}\right)}
\end{aligned}
\end{equation}
The analytical solution for the displacement is given by:
\begin{equation}
u(z,t)=\frac{3P_{A}}{3K+4G}\Bigg[(h-z)+\gamma \alpha \sum \limits_{n=1,3}^\infty  E_{1-\beta,1} \left(-{\dfrac{n^{2}\pi^{2} \bar{\lambda} t^{1-\beta}}{4h^{2}}}\right) \dfrac{8h}{(n \pi)^2}\bigg\{(-1)^{\dfrac{n-1}{2}} \sin \frac{n \pi z}{2h}-1 \bigg\} \Bigg]
\label{eq:DispSoln}
\end{equation}
The detailed derivation is given in Appendix A. \\
The expression for the flux is found by substituting the solution for the pore pressure (c.f. \eref{eq:PorepSoln}) in the flux definition (c.f. \eref{eq:FluidFlux}) and solving the fractional derivative as follows:
\begin{equation}
J_p = \lambda_\beta P_A \gamma t^{-\beta}\sum \limits_{n=1,3}^\infty E_{1-\beta,1-\beta}\Bigg(\dfrac{-n^2\pi^2\bar{ \lambda}t^{1-\beta}}{4h^2} \Bigg)\dfrac{2}{h} \sin\left(\dfrac{n\pi z}{2h}\right)
\label{eq:fluxfullform}
\end{equation}
The weight of the sample over time $W(t)$ is related to the initial weight of the sample $W_{0}$ and the fluid flux $J_p$ as follows:
\begin{equation}
W(t)= W_{0} - \int\limits_{0}^{t}J_p \cdot w_s ~\boldsymbol{A}~\mathrm{d}t  \\
\label{eq:WtLoss}
\end{equation}
where $w_s$ is the specific weight of the fluid, and $\boldsymbol{A}$ is the cross-section over which the fluid flows out. Substituting  \eref{eq:fluxfullform} in \eref{eq:WtLoss} and evaluating the integral analytically, we obtain the following relation:
\begin{equation}
W(t)= W_{0}- P_{A} A w_s\gamma \frac{2}{h}\sum \limits_{n=1,3}^\infty  \left[\lambda_{\beta}t^{1-\beta} E_{1-\beta,2-\beta} \left(-{\frac{\pi^{2} \bar{\lambda} t^{1-\beta}}{4h^{2}}}\right)\right]    \\
\label{eq:21}
\end{equation}
For detailed derivation, refer to Appendix B.

\begin{figure}[htp!]
\begin{subfigure}{.5\textwidth}
\centering
\includegraphics[width=\linewidth]{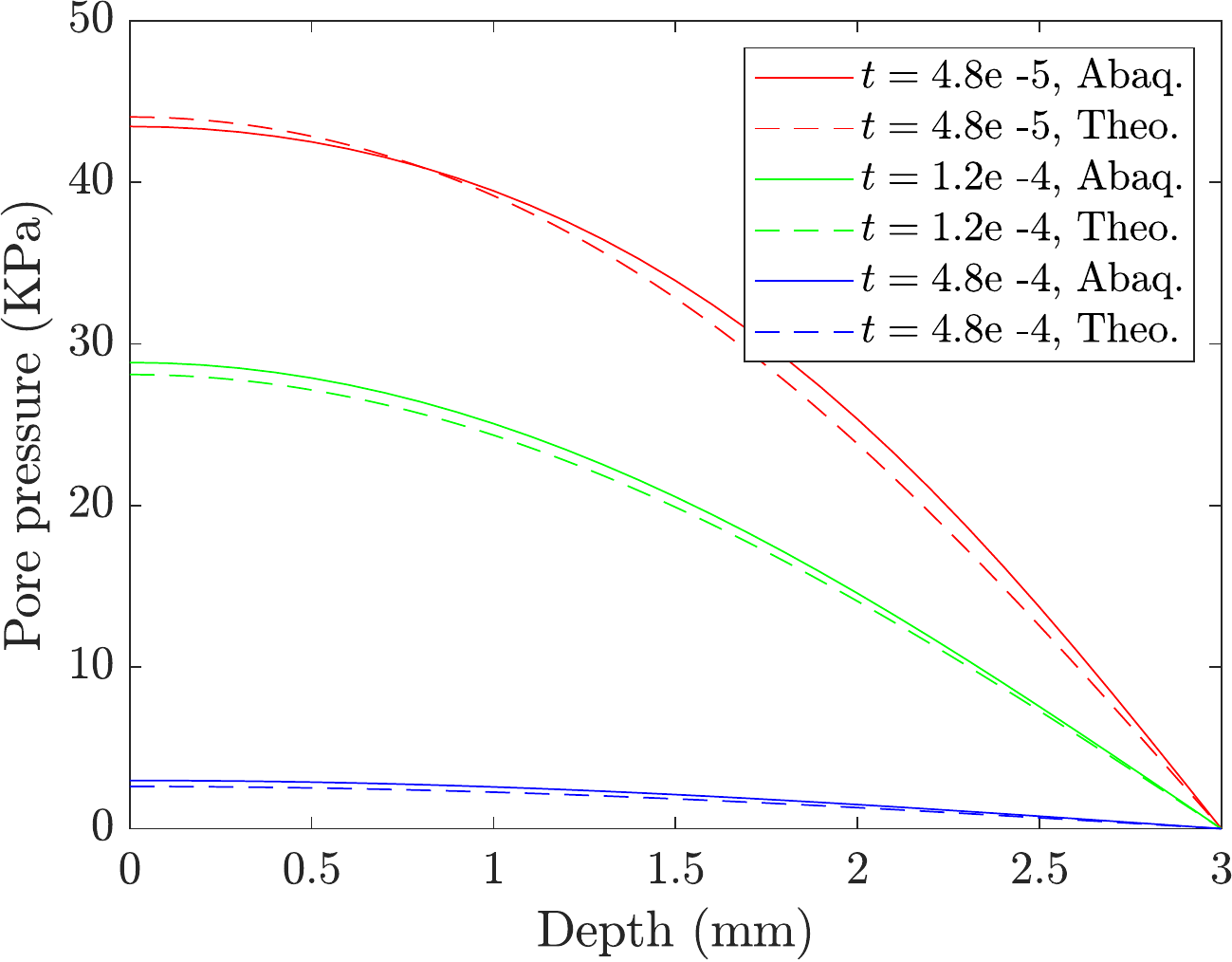}
\caption{\label{fig:pvsd_AbqTh_comp}}
\end{subfigure}
\begin{subfigure}{.5\textwidth}
\centering
\includegraphics[width=\linewidth]{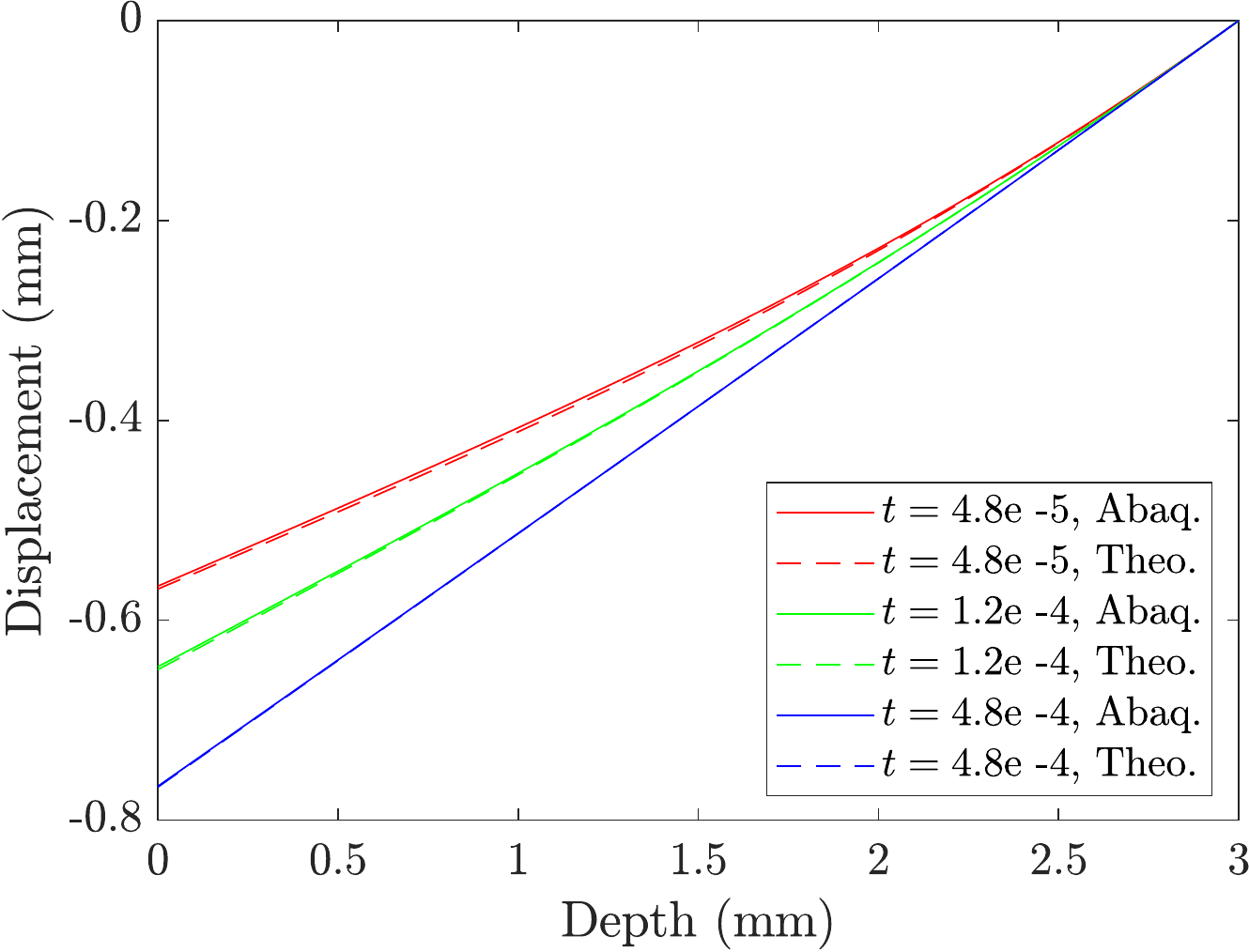}
\caption{\label{fig:dvsd_AbqTh_comp}}
\end{subfigure}\\
\begin{subfigure}{.5\textwidth}
\centering
\includegraphics[width=\linewidth]{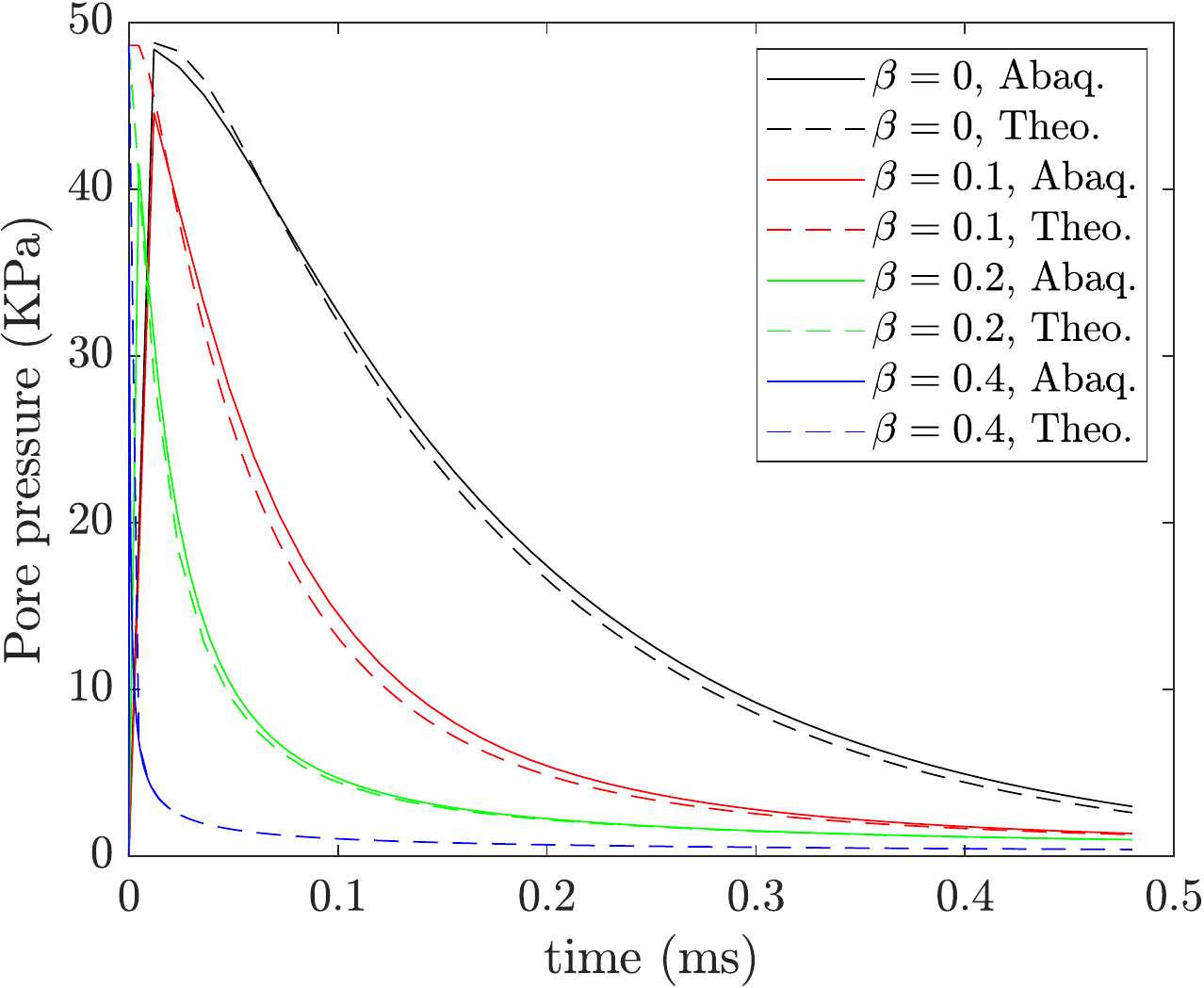}
\caption{\label{fig:pvst_AbqTh_comp}}
\end{subfigure}%
\begin{subfigure}{.5\textwidth}
\centering
\includegraphics[width=\linewidth]{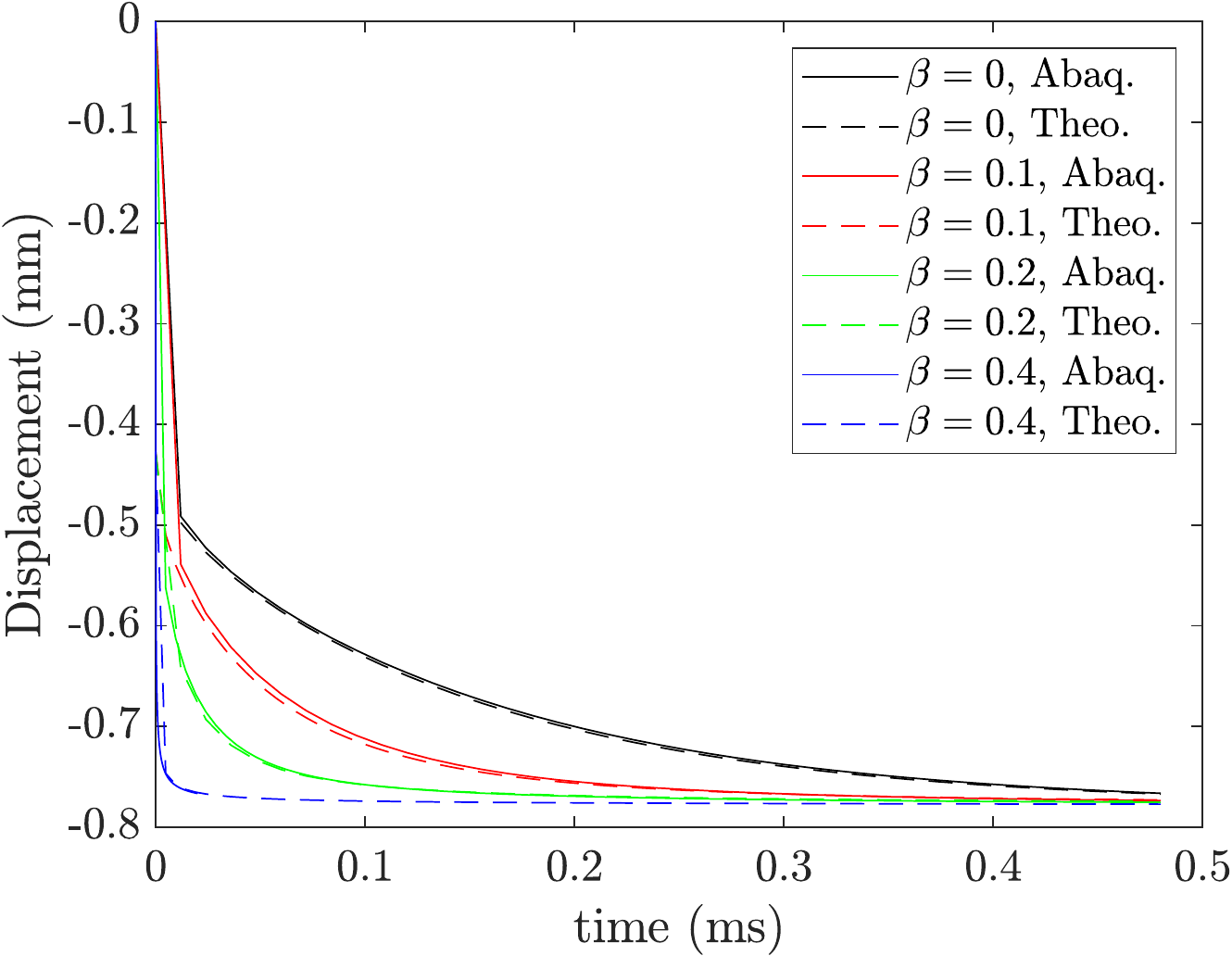}
\caption{\label{fig:dvst_AbqTh_comp}}
\end{subfigure}
\caption{Comparison of numerical results with theoretical results: (a-b) variation of pore pressure and displacement with depth at different time steps and $\beta=$ 0 and (c-d) influence of $\beta$ on the pore pressure and displacement as a function of time for a point on the top surface.}
\label{fig:AbqTh_comp}
\end{figure}

\paragraph{Finite Element implementation and verification}
Biot's model with fractional Darcy's law is implemented in Abaqus using UMATHT and UMAT subroutines using the method proposed by Barrera~\citep{Barrera2021}.In ~\citep{Barrera2021} only the pore pressure field is discussed, here for the first time we show that both the pore pressure and the displacement fields agree with the analytical solution. Using the similarity between the governing equations of thermoelasticity and poroelasticity~\citep{Barrera2021}, UMATHT of Abaqus, which is generally meant for thermoelasticity, is used for the case of poroelasticity. In Abaqus, Temperature can be treated as Pore pressure, and the heat flux can be treated as heat flux. UMAT is used for transferring the stress data to the UMATHT for modelling the coupling behaviour. UVARM can be used for transferring stress values. However, due to the order of implementation in Abaqus, which is elastic model-UMATHT-UVARM, the stress values from UVARM reach the UMATHT in the next iteration, which introduces errors in the computation. 
To avoid this error, UMAT is used, which allows for sharing stress values in the same increment. UMAT is written for a linear elastic model. The coefficient coupling the strain and pore pressure $\alpha_T = \alpha/3K$ is given as the thermal expansion coefficient.
\begin{table}[htpb]
\caption{Material parameters}
\centering
\begin{tabular}{ l  l }
\hline
Bulk modulus, $K$ & 1.67$\times$10$^{5}$Pa\\  
Shear modulus, $G$ & 7.69$\times$ 10$^4$ Pa ($E= $     0.2$\times$10$^{6}$Pa,  $\nu=$ 0.3)\\
Skempton coefficient, $B$ & 0.88 \\  
Biot coefficient, $\alpha$ & 0.65 \\  
Permeability, $\lambda_\beta$ & 8.33$\times$10$^{-8}$m$^4$/Ns$^{1-\beta}$\\ 
Fractional Order, $\beta$ & 0, 0.1, 0.2, 0.4\\ \hline
\end{tabular}
\label{T:Material properties compression case}
\end{table}
Abaqus model is validated by comparing with theoretical solutions of a confined compression test with material properties given in \tref{T:Material properties compression case}~\citep{Barrera2021}. For this, a cylindrical computational model with height $h=$ 3mm and diameter $d=$ 3mm is considered. The domain is discretized with 8-noded trilinear hexahedral elements with 4 degrees of freedom per node (C3D8T elements). The displacements in all three directions are restrained on the bottom face and the displacements on the side faces are restrained laterally. The pore pressure is zero at the bottom and an instantaneous load, $P_A=$ 0.07 MPa is applied on the top surface. Based on a systematic study, a total of 28,110 8-noded hexahedral elements, were found to be adequate to model the poroelastic response. The displacements and the pore pressure obtained from the numerical simulation are compared against the theoretical results. \frefs{fig:pvsd_AbqTh_comp}-\ref{fig:dvsd_AbqTh_comp} shows the comparison of the pore pressure and displacement as a function for depth for different time steps for $\beta=$ 0 between theoretical solution and numerical computation. The influence of $\beta$ on the pore pressure and the displacement is shown in \frefs{fig:pvst_AbqTh_comp}-\ref{fig:dvst_AbqTh_comp}. It is seen from \fref{fig:AbqTh_comp} that there is a very good agreement with the numerical result and theoretical solution and that the fractional order has a strong influence on the pore pressure and the displacement. Large jumps in the pore pressure and displacement in the initial time in Figures \ref{fig:pvst_AbqTh_comp}, \ref{fig:dvst_AbqTh_comp} are due to the numerical approximation of step function with ramp in the first time increment.


\section{Results}
\subsection{Material parameters: anomalous permeability, order of derivative and aggregate modulus}
\begin{figure}[H]
\centering
\begin{subfigure}{.5\textwidth}
\centering
\includegraphics[width=\linewidth]{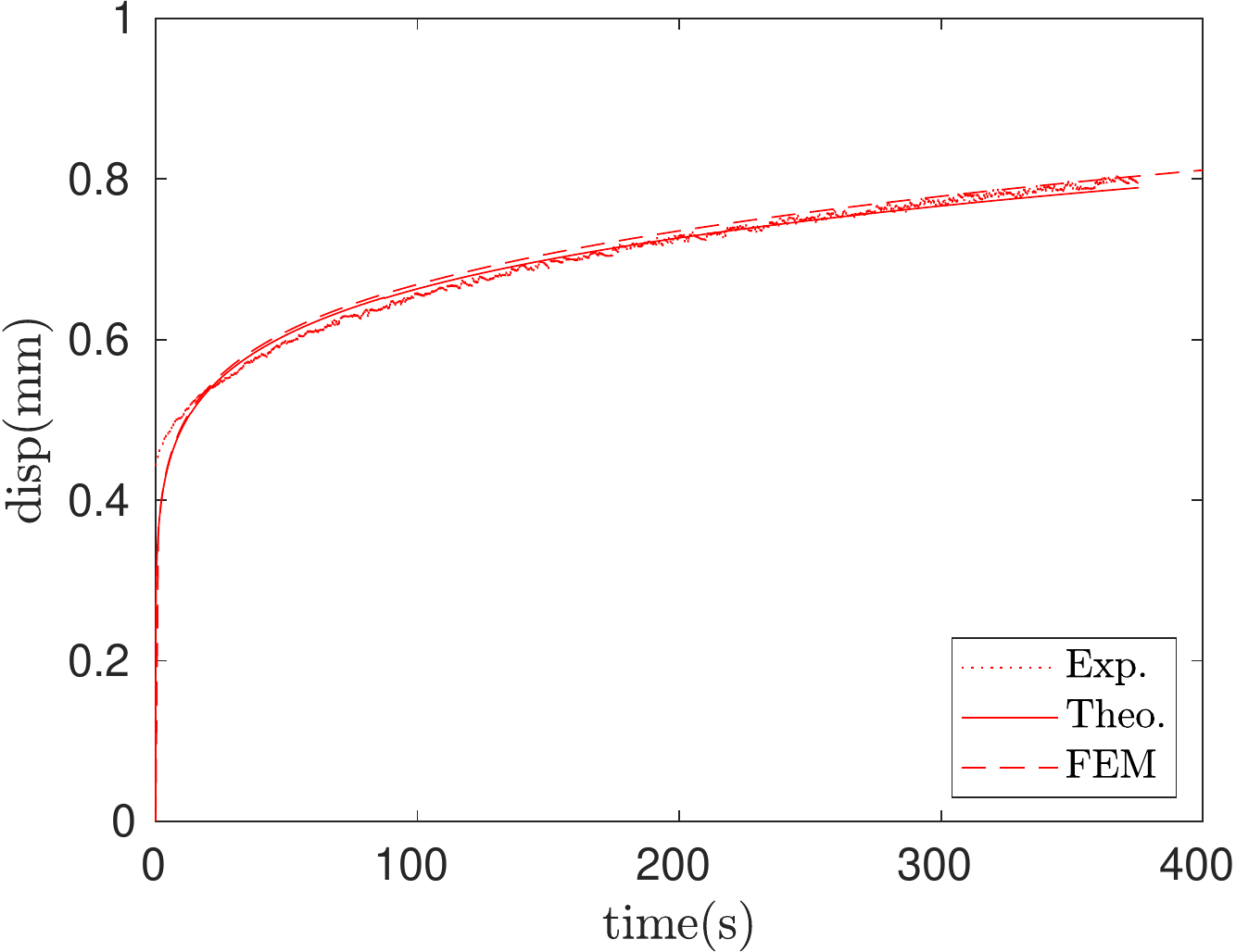}
\caption{TK11 Body Circumferential}
\label{TK11BC_dvst_AbqExpTh}
\end{subfigure}%
\begin{subfigure}{.5\textwidth}
\centering
\includegraphics[width=\linewidth]{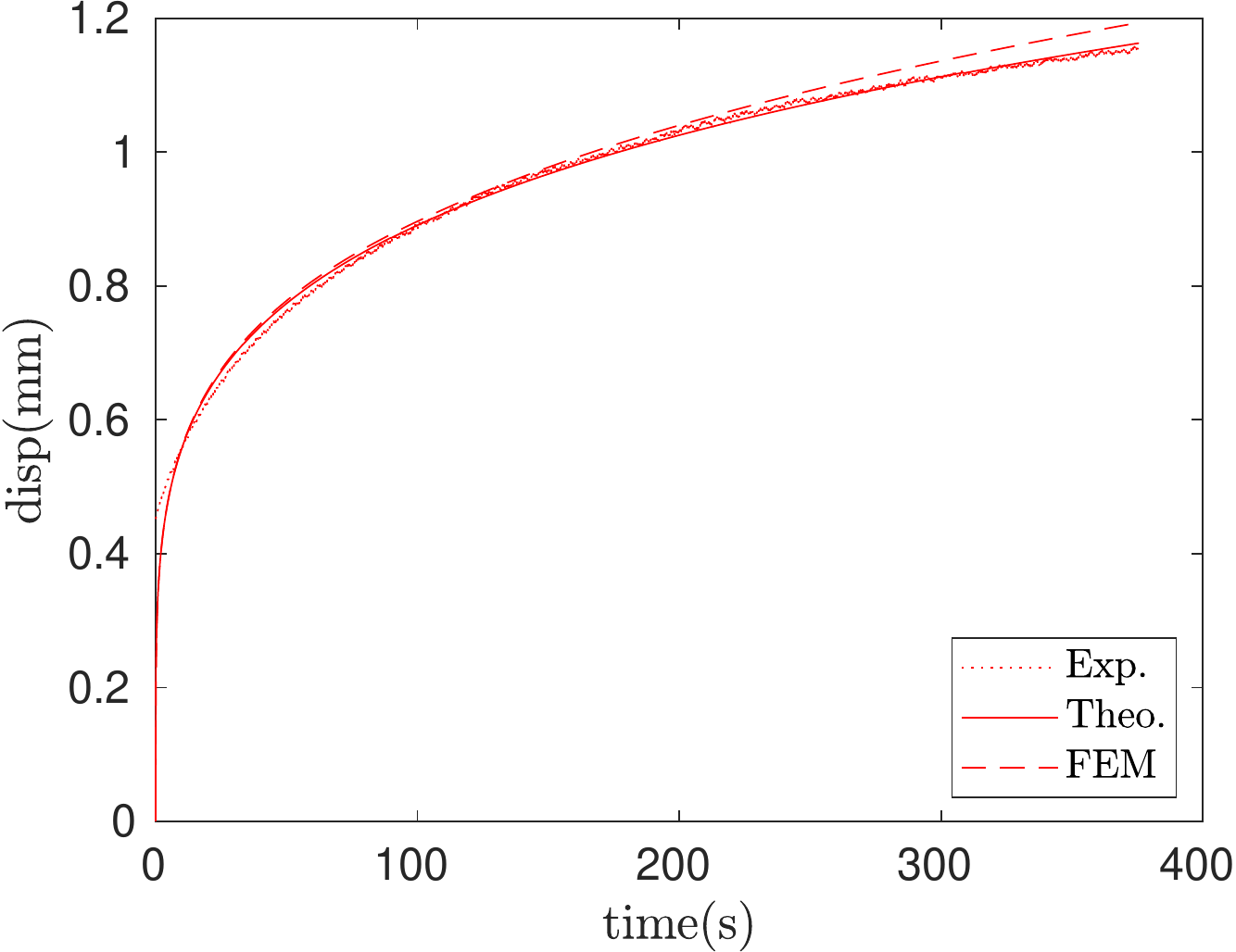}
\caption{TK11 Body Radial}
\label{TK11BR_dvst_AbqExpTh.pdf}
\end{subfigure}\\%
\begin{subfigure}{.5\textwidth}
\centering
\includegraphics[width=\linewidth]{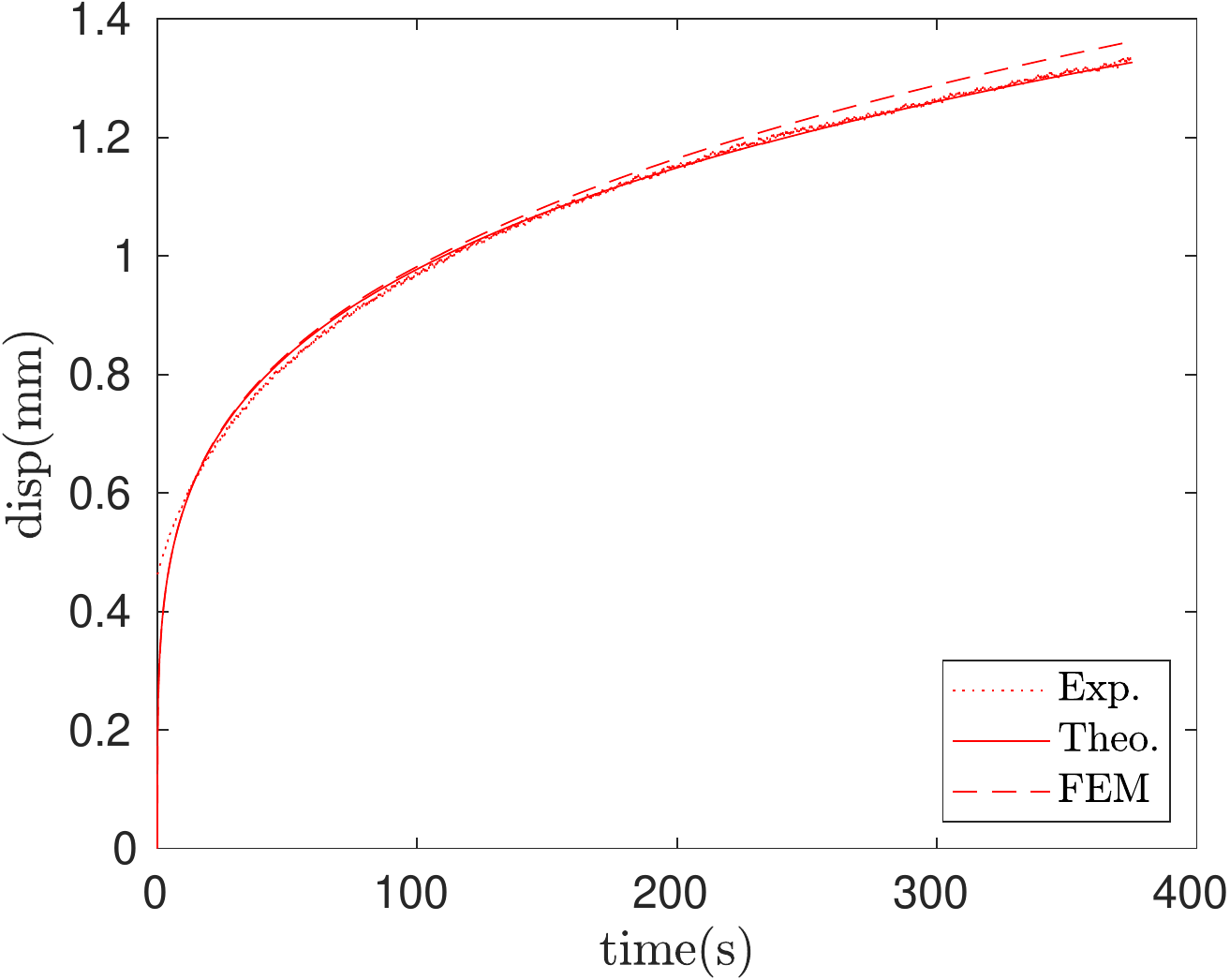}
\caption{TK11 Body Vertical}
\label{TK11BV_dvst_AbqExpTh.pdf}
\end{subfigure}%
\begin{subfigure}{.5\textwidth}
\centering
\includegraphics[width=\linewidth]{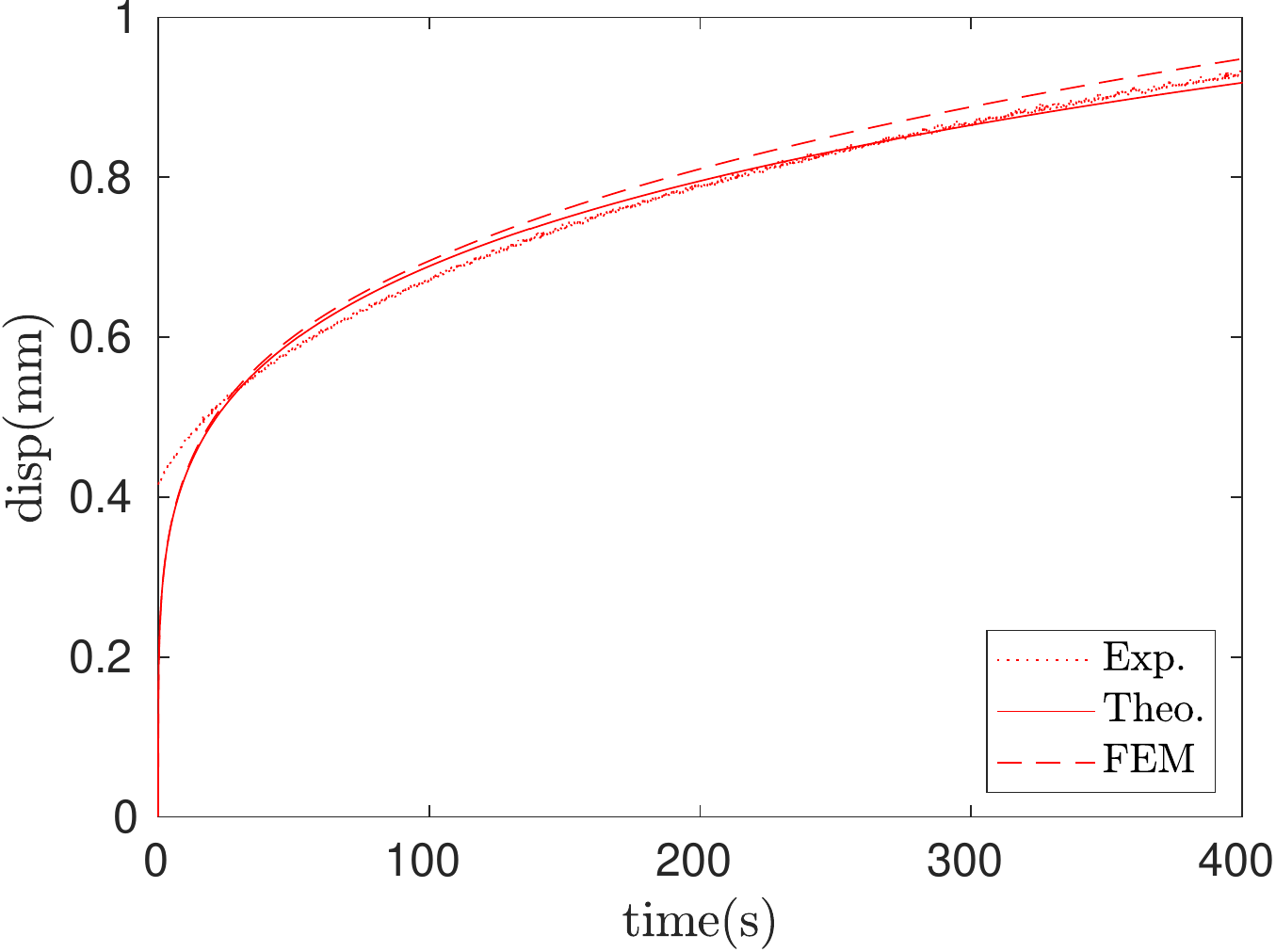}
\caption{TK16 Anteriorhorn Radial}
\label{TK16AR_dvst_AbqExpTh.pdf}
\end{subfigure}\\%
\caption{Displacement confined creep data fitting with fractional poroelasticity and modelling using FEM, (a) Sample taken from the body region of TK11 sample in the circumferential direction, (b) Sample taken from the body region of TK11 sample in the radial direction, (c) Sample taken from the body region of TK11 sample in the vertical direction, (d) Sample taken from the anterior horn region of TK16 sample in the radial direction.  }
\label{fig:CreepDispCompar}
\end{figure}

The proposed fractional poroelastic model is used to characterize the response of the meniscus. To this, we use the expression in \eref{eq:DispSoln} for fitting the experimental data. The meniscus is assumed to be incompressible material, due to which $B = 1, \alpha = 1,$ and $K_u \rightarrow\infty$ reducing \eref{eq:DispSoln} to:
\begin{equation}
u(z,t)=\dfrac{P_{A}}{M}\Bigg[(h-z)+ \sum \limits_{n=1,3}^\infty  E_{1-\beta,1} \left(-{\dfrac{n^{2}\pi^{2} \lambda_\beta M t^{1-\beta}}{4h^{2}}}\right) \dfrac{8h}{(n \pi)^2}\bigg((-1)^{\dfrac{n-1}{2}}\sin \dfrac{n \pi z}{2h}-1 \bigg) \Bigg]
\label{eq:DispIncom}
\end{equation}
where, $M=(3K + 4G)/3$, is the aggregate modulus. \eref{eq:DispIncom} is used to find the material properties of meniscus tissue by fitting with the creep displacement vs time data. It is noticed that this model requires only three material properties, viz., the aggregate modulus ($M$), the fractional order ($\beta$) and the permeability ($\lambda_\beta$). On using $\beta=$ 0, this theory converts into the classical Biot's theory. The experiments were displacement controlled, i.e., a constant piston velocity was used until the load reached the required value. So, ramp loading is neglected for fitting and approximated as step loading. Displacements measured in the experiments were fitted with the analytical solution \eref{eq:DispIncom} at $z=0$, i.e.,
\begin{equation}
u = \frac{P_A}{M} \bigg[h - \sum_{n=1,3}^\infty E_{1-\beta,1}\bigg(\frac{-n^2\pi^2 \lambda_\beta Mt^{1-\beta}}{4h^2}\bigg)\frac{8h}{(n\pi)^2}\bigg]  
\end{equation}
Fitting is done in MATLAB using the inbuilt function `fminsearch' with material properties as variables and the RMS error between the analytical and the experimental data as the function to be minimized. Fitting plots for four samples are shown in \fref{fig:CreepDispCompar}.
\tref{T:MaterialProperties} shows the material properties of the meniscus obtained from the experimental results for samples in both body and anterior horn for three different directions, viz., circumferential, vertical and radial. The last two characters of the sample name indicate the location and orientation information, for example, BC represents Body Circular.  The other notations employed are A: Anterior horn, P: Posterior horn, R: Radial, and V: vertical directions. 
\begin{figure}[H]
\begin{subfigure}{\textwidth}
\centering
\includegraphics[width=0.5\linewidth]{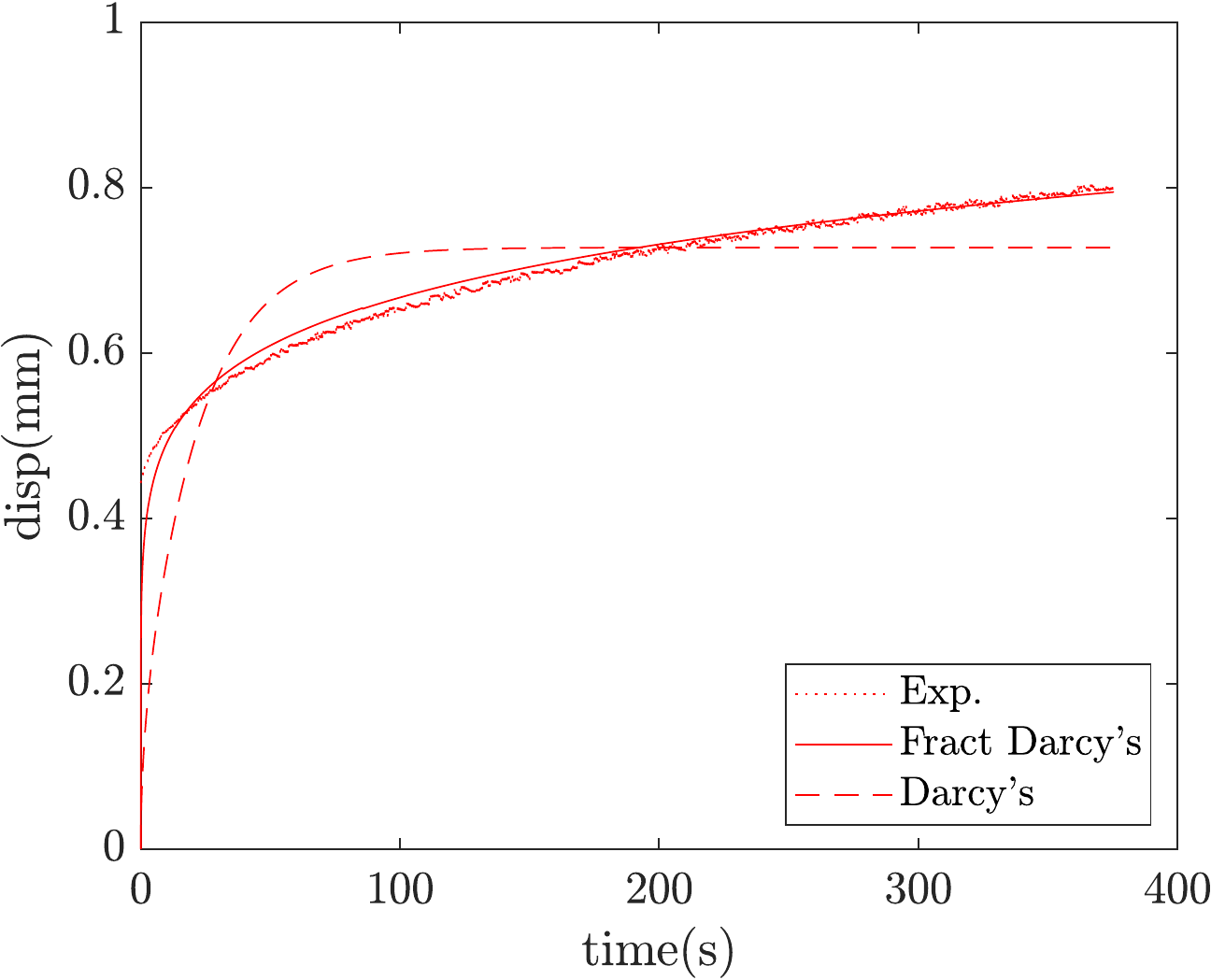}
\caption{}
\label{}
\end{subfigure}\\%
\begin{subfigure}{0.5\textwidth}
\centering
\includegraphics[width=\linewidth]{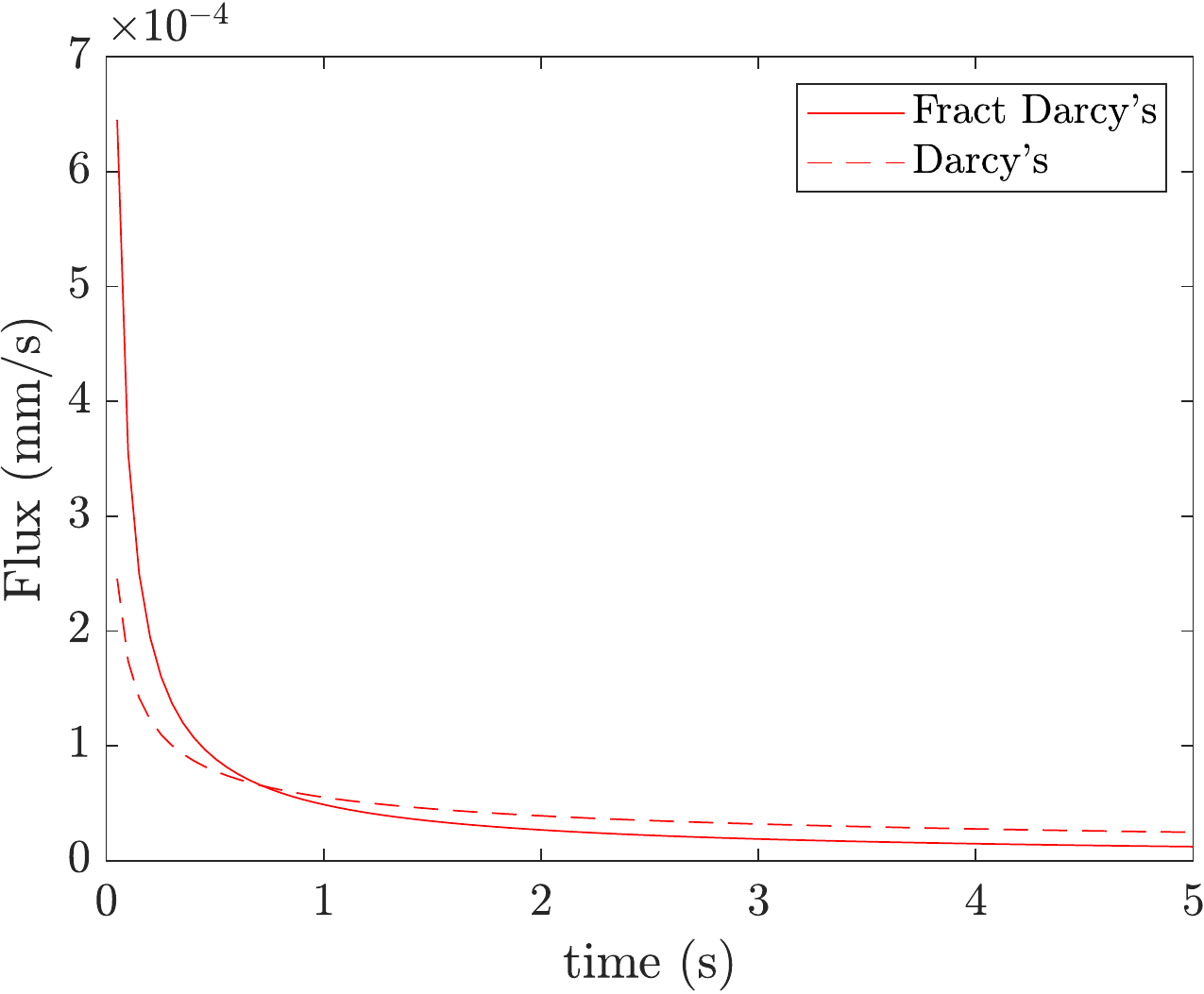}
\caption{}
\label{}
\end{subfigure}%
\begin{subfigure}{.5\textwidth}
\centering
\includegraphics[width=\linewidth]{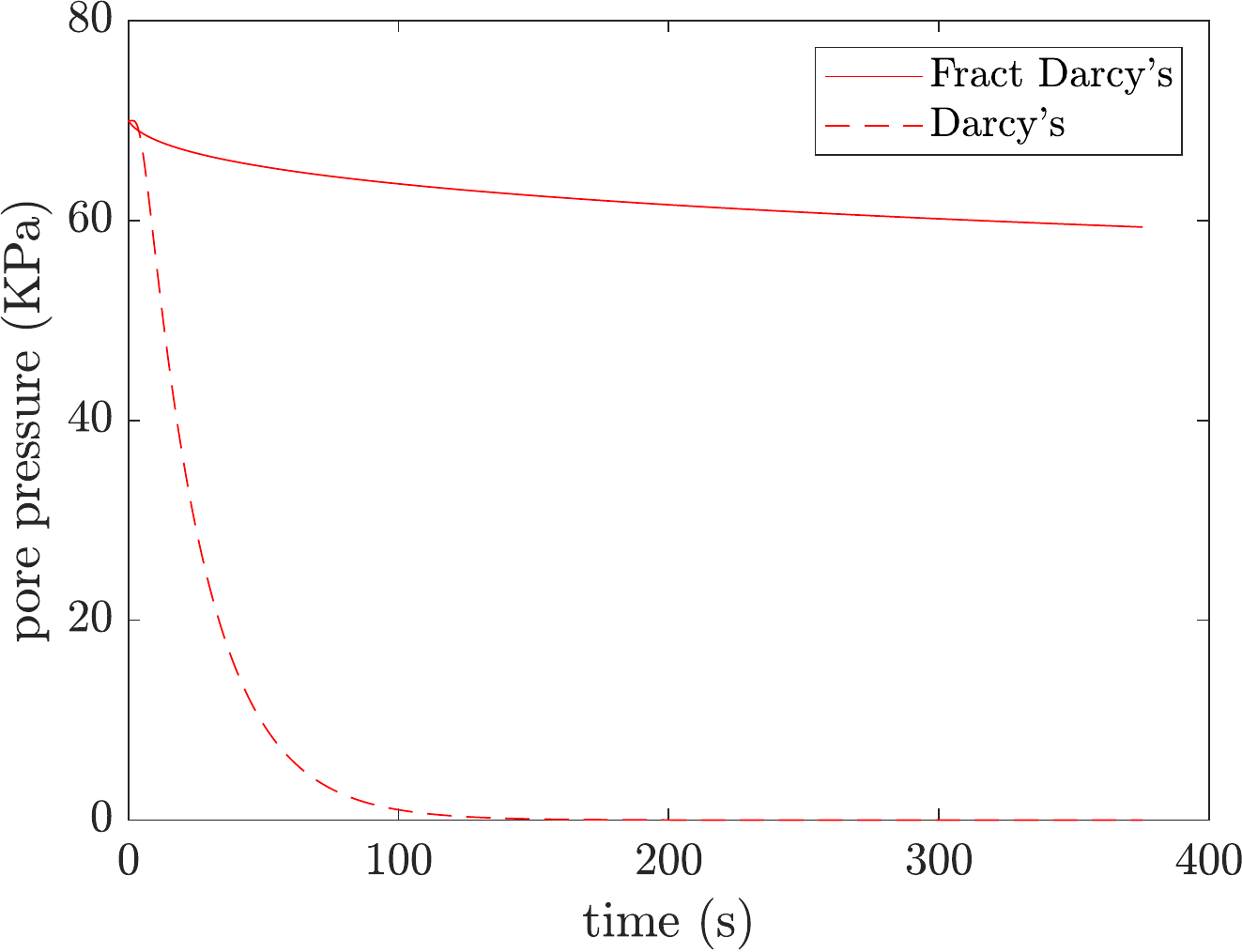}
\caption{}
\label{}
\end{subfigure}%
\caption{Comparison between Biot's theory with classical and fractional Darcy's law to fit creep displacement experiment data.}
\label{fig:ClassicvsFrac}
\end{figure}

To understand the fractional order's influence, we compared the classical Biot's poroelastic model and the proposed one (c.f. \sref{sec:governing eqns}). The material properties obtained for one sample TK11BC through a fitting for classical Biot's model\footnote{Refer Appendix  \ref{A:MatProps_classic} for material properties corresponding to the classical Biot's theory for all the samples.}  are $M =$ 3.56$\times$10$^{-1}$ MPa \& $\lambda =$ 6.95$\times$10$^{-13}$ m$^4$/Ns and whilst that with the fractional framework are: $M =$ 1.27$\times$10$^{-1}$ MPa, $\beta =$ 0.73 \& $\lambda_{\beta} =$ 2.95$\times$10$^{-12}$ m$^4$/Ns$^{1-\beta}$. The  RMS error using Classical Biot's theory is 6.53$\times$10$^{-5}$, and the RMS error while using  Biot's theory with fractional Darcy's law is 1.42$\times$10$^{-5}$. \fref{fig:ClassicvsFrac} shows a comparison of the displacement, flux out and pore pressure as a function of time for classical Biot's theory and fractional model with experiments. It is opined that the classical Biot's model does not fit the experimental results well and reemphasizes a need for a fractional poroelastic framework.
        
\begin{table}[!htp]\centering
\caption{Meniscus material parameters obtained from fitting with creep data}
\label{Menscus_Param }
\scriptsize
\renewcommand{\arraystretch}{1.2}
\begin{tabular}{llcccrc}
\toprule
\textbf{S No.} &\textbf{Sample} & \textbf{h}&\textbf{M}$\times$10$^5$ &\textbf{$\boldsymbol\beta$} &\textbf{$\boldsymbol{\lambda_\beta}$}$\times$10$^{-12}$ & \textbf{RMS error} \\
 & & (mm) & (Pa) & & (m$^2$/Pa.s$^{1-\beta})$ & ($\times$10$^{-5}$)\\ 
\midrule
1	    & TK11BC    & 3.7	&   1.27	& 0.73	& 2.95	& 1.42\\
2	    & TK11BR    & 3.3	&   0.71	& 0.59	& 1.48	& 1.61\\
3	    & TK11BV	& 4.1	&   0.83	& 0.53	& 1.56	& 1.74\\
4	    & TK16BC1	& 3.2	&   0.76	& 0.7	& 2.13	& 1.39\\
5	    & TK16BR2	& 2.9	&   0.82	& 0.59	& 0.89	& 1.30\\
6	    & TK16BV	& 3.2	&   0.62	& 0.59	& 0.85	& 1.48\\
7	    & TK16AC    & 3.2	&   0.60	& 0.63	& 0.97	& 1.97\\
8	    & TK16AR1	& 3.33	&   0.85	& 0.58	& 1.01	& 2.02\\
9	    & TK16AV	& 3.33	&   0.57	& 0.5	& 0.37	& 1.63\\
10	    & TK16BC2	& 3.33	&   0.90	& 0.76	& 4.28	& 1.53\\
11	    & TK16BR1	& 2.9	&   0.55	& 0.64	& 0.91	& 1.83\\
12	    & TK16BV2	& 3.13	&   0.61	& 0.63	& 0.93	& 1.56\\
13	    & TK16PV	& 2.67	&   1.41	& 0.75	& 1.99	& 1.40\\
14	    & TK17BC	& 3.2	&   0.10	& 0.73	& 1.03	& 2.58\\
15	    & TK17BR	& 2.2	&   0.41	& 0.72	& 1.36	& 1.55\\
16	    & TK17BV	& 2.3	&   0.86	& 0.54	& 0.90	& 1.31\\
17	    & TK17AC	& 2.73	&   1.87	& 0.79	& 4.96	& 1.49\\
18	    & TK17AR	& 2.7	&   1.94	& 0.76	& 3.05	& 1.44\\
19	    & TK18BC	& 4.6	&   1.34	& 0.74	& 141.03	& 3.81\\
20	    & TK18BR	& 2.9	&   0.66	& 0.72	& 1.90	& 1.78\\
21	    & TK18BV	& 2.8	&   0.68	& 0.61	& 1.07	& 1.63\\
22	    & TK22BR	& 3.17	&   2.54	& 0.49	& 85.2	& 9.90\\
23	    & TK22AC	& 3.67	&   0.65	& 0.62	& 1.40	& 2.68\\
24	    & TK22AR	& 3.33	&   1.60	& 0.41	& 2.24	& 4.69\\
25	    & TK22AV	& 2.93	&   0.76	& 0.51	& 0.35	& 1.42\\
26	    & TK36BC	& 3.37	&   0.74	& 0.73	& 4.40	& 2.11\\
27	    & TK37BC	& 4.6	&   0.90	& 0.74	& 6.63	& 3.53\\
 28	    & TK37BR	& 3.5	&   0.67	& 0.59	& 0.96	& 1.93\\
29	    & TK37BV	& 3.2	&   0.49	& 0.53	& 0.55	& 1.33\\
\bottomrule
\end{tabular}
\label{T:MaterialProperties}
\end{table}

\begin{figure}[H]
\centering
\begin{subfigure}{.3\textwidth}
\centering
\includegraphics[width=\linewidth]{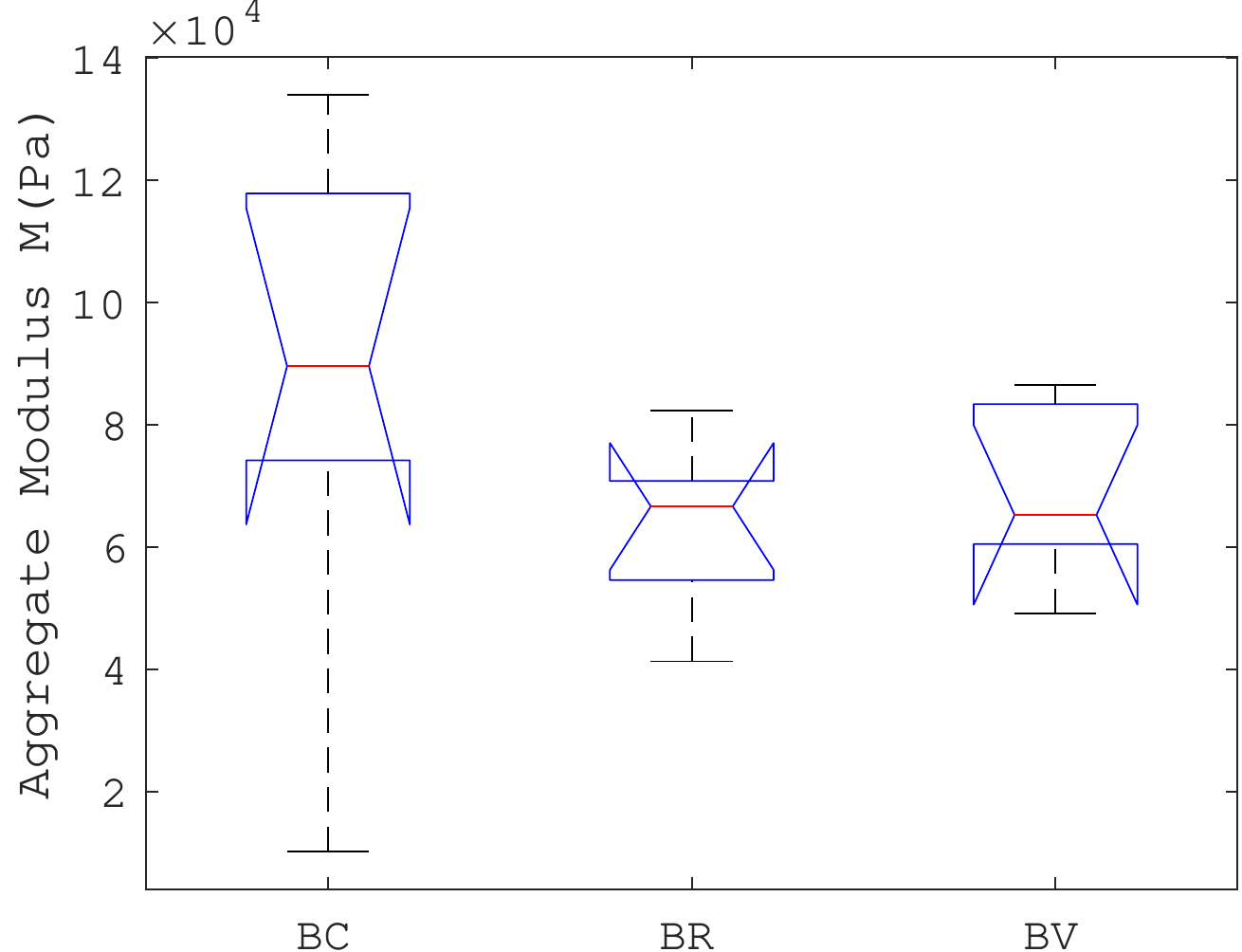}
\caption{}
\label{subfig:M_Body_Anova}
\end{subfigure}  %
\begin{subfigure}{.3\textwidth}
\centering
\includegraphics[width=\linewidth]{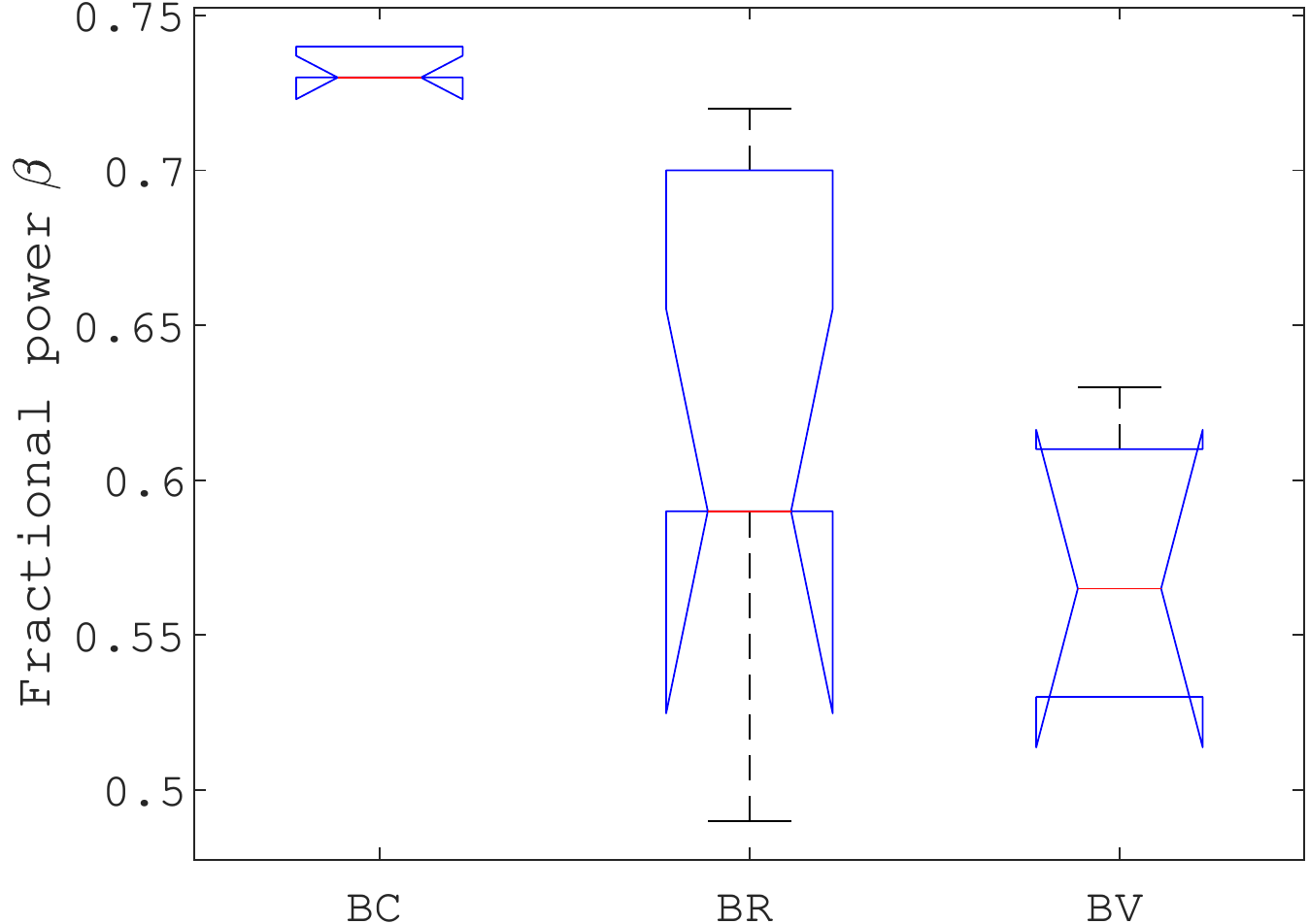}
\caption{}
\label{subfig:Beta_Body_Anova}
\end{subfigure}  %
\begin{subfigure}{.3\textwidth}
\centering
\includegraphics[width=\linewidth]{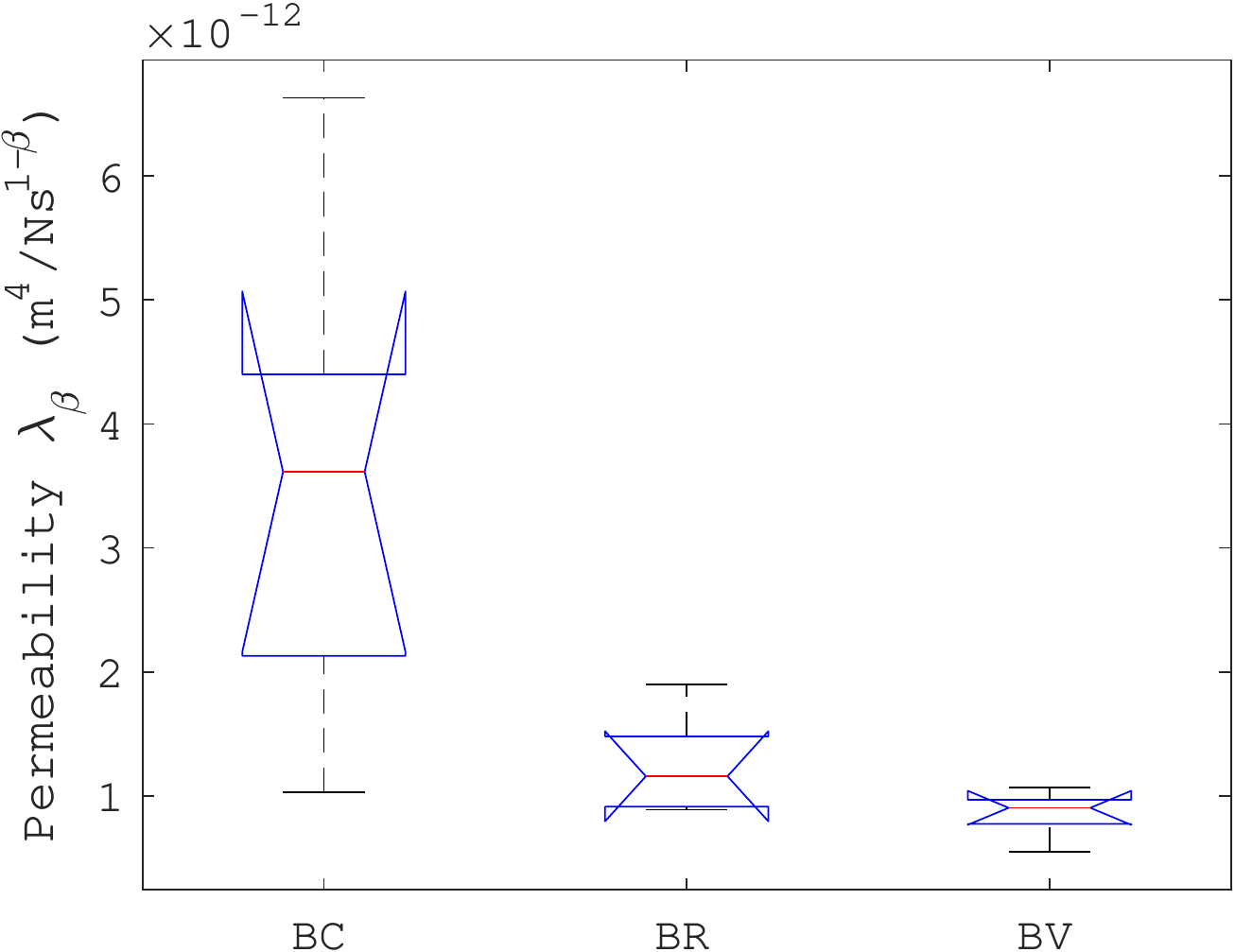}
\caption{}
\label{subfig:Lambda_Body_Anova}
\end{subfigure}\\%
\caption{Material properties comparison in different orientations in the body region(BC - Body circumferential, BR - Body Radial, BV - Body Vertical). (a) Aggregate Modulus (M), (b) Fractional power $\beta$, (c) Permeability $\lambda_\beta$}
\label{fig:Body_Anova}
\end{figure}

    
Further, to check for anisotropy in the meniscus, the ANOVA test is performed on the body region, properties with circular, radial and vertical directions as categories and properties are compared in different directions using MATLAB's default function `anova1'. For a few samples, properties were abnormal compared to the remaining samples of the same category. It could be due to the sample, experiment or fitting issue. These samples were detected as outliers by the Matlab 'anova1' function and were not considered for comparison. The mean and the standard deviation of the groups of the samples without considering the outliers are given in \tref{Meniscus param sample} and the results for the ANOVA test are shown in \fref{fig:Body_Anova}. In \fref{subfig:M_Body_Anova}, the aggregate modulus is compared. The red line shows the sample's median, and the black dashed line limit gives the sample's range. If the notched blue box of different samples does not overlap, it is concluded that the true medians are different with 95\% confidence. If the probability against the null hypothesis ($p$-value) is less than 0.05, it is concluded that the mean of the categories is different with 95\% confidence. For aggregate modulus, the obtained $p$-value is 0.329, indicating that
the aggregate modulus is not different in different directions and could possibly be treated as isotropic. The obtained $p$-value for fractional order $\beta$ and permeability $\lambda_\beta$ are 0.001 and 0.0038, respectively. From \frefs{subfig:Beta_Body_Anova}-\ref{subfig:Lambda_Body_Anova}, it is inferred that the fractional order and permeability are relatively higher in circular directions when compared to the other two directions, viz., vertical and radial. Further, it can be considered similar in vertical and radial directions.

\begin{table}[!htp]\centering
\caption{Meniscus parameters variation across sample}\label{Meniscus param sample}
\scriptsize
\begin{tabular}{lrrrrrrr}\toprule
&\multicolumn{2}{c}{$\boldsymbol{M}\times$10$^5$} &\multicolumn{2}{c}{$\boldsymbol\beta$} &\multicolumn{2}{c}{$\boldsymbol{\lambda_\beta}\times$10$^{-12}$} \\
&\multicolumn{2}{c}{(Pa)} &\multicolumn{2}{c}{} &\multicolumn{2}{c}{(m$^2$/Pa.s$^{1-\beta}$)}\\
\cmidrule{2-7}
Part &Mean &SD &Mean &SD &Mean &SD \\\midrule
Body Cir &0.75 &0.49 &0.73 & 0.00 &3.75 &2.36 \\
Body Rad &0.64 &0.14 &0.64 &0.06 &1.25 &0.40 \\
Body Ver &0.65 &0.14 &0.58 &0.04 &0.86 &0.19 \\
Anthorn Cir &1.04 &0.72 &0.68 & 0.10 &2.44 & 2.19 \\
Anthorn Rad &1.47 &0.56 &0.58 &0.17 &2.10 &1.03 \\
Anthorn Ver &0.66 &0.14 &0.51 &0.01 &0.36 &0.02 \\
\bottomrule
\end{tabular}
\end{table}
           
\begin{figure}[H]
\begin{subfigure}{.5\textwidth}
\centering
\includegraphics[width=\linewidth]{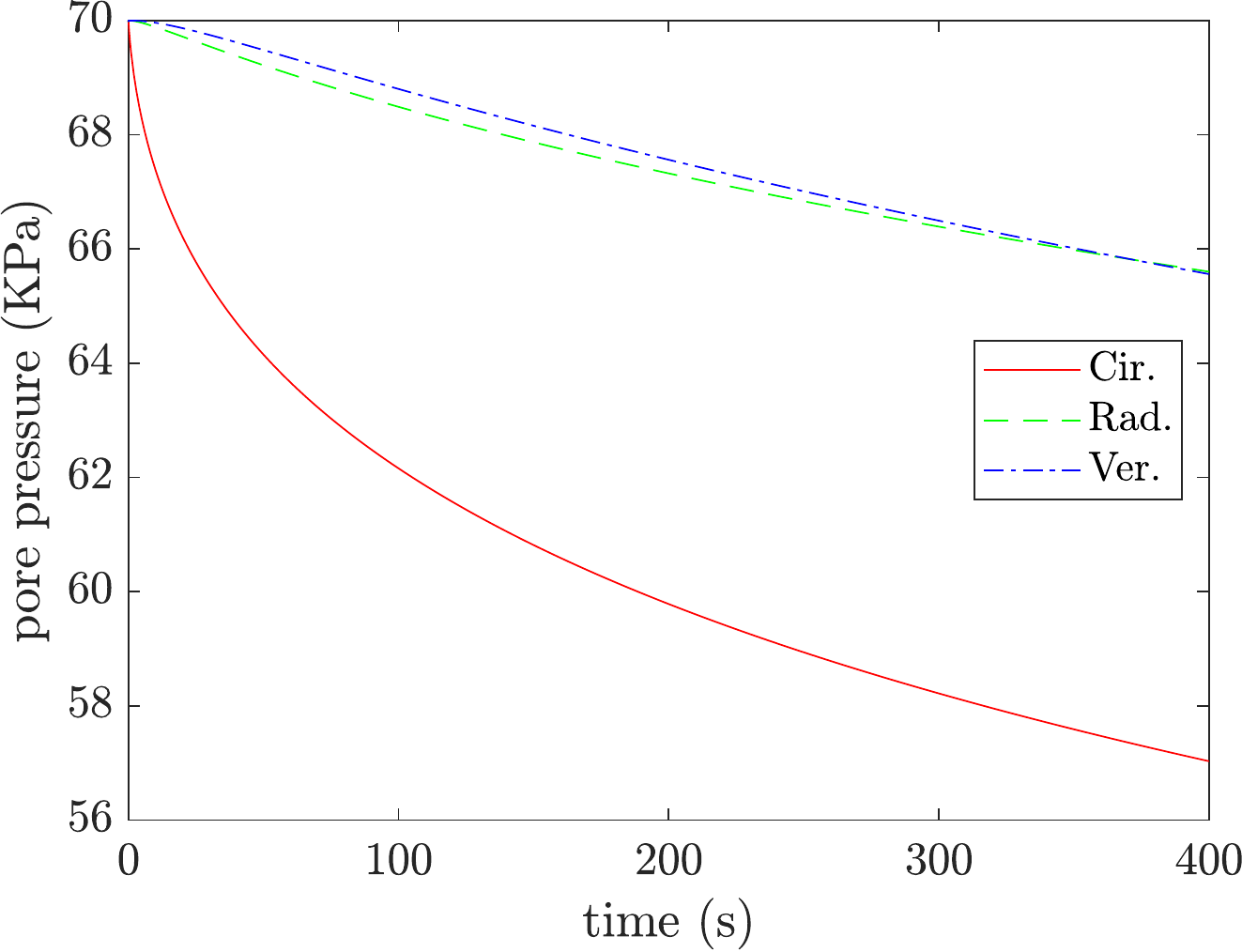}
\caption{}
\label{sf:porep_directional_compar}
\end{subfigure}
\begin{subfigure}{.5\textwidth}
\centering
\includegraphics[width=\linewidth]{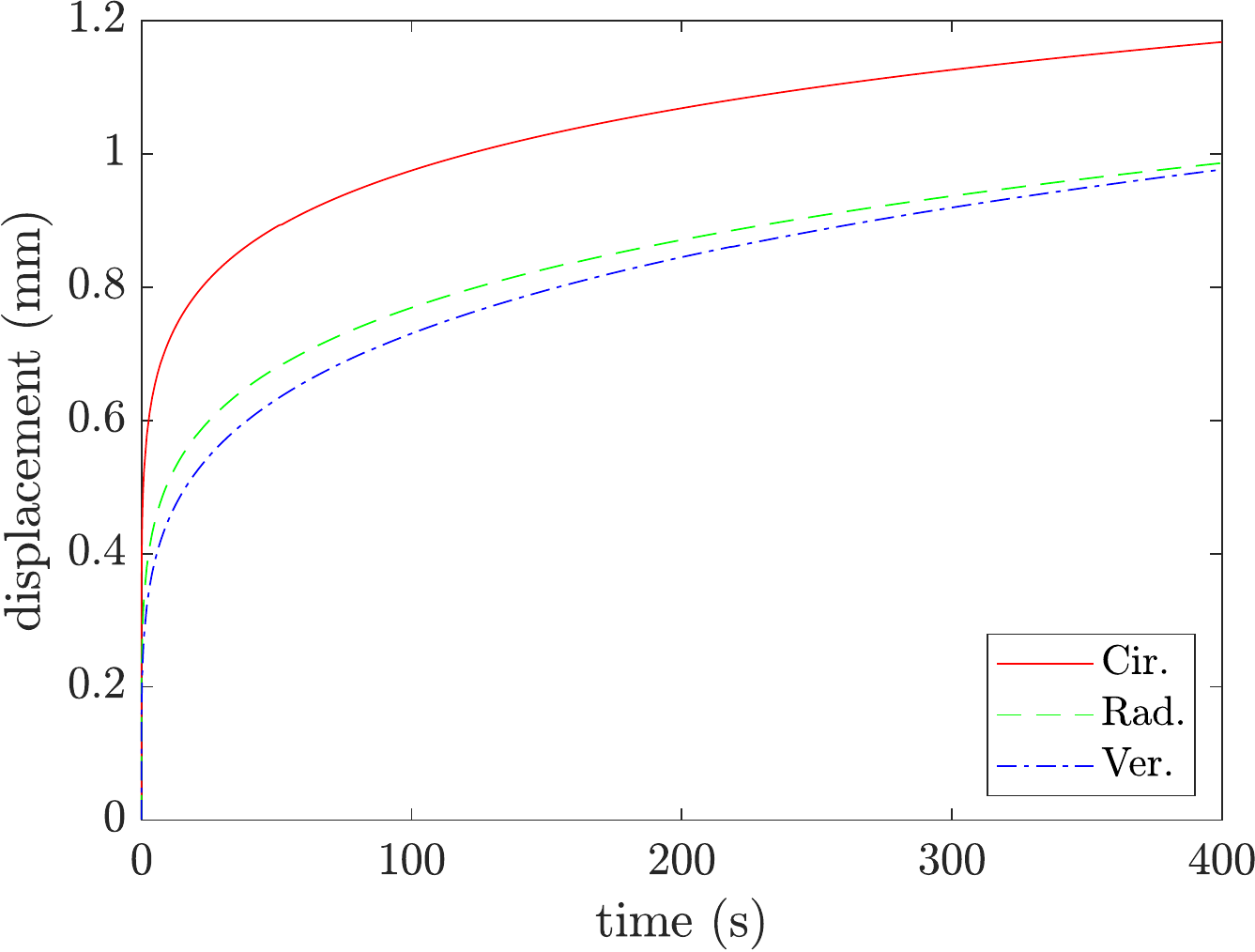}
\caption{}
\label{sf:disp_directional_compar}
\end{subfigure}\\%
\begin{subfigure}{.5\textwidth}
\centering
\includegraphics[width=\linewidth]{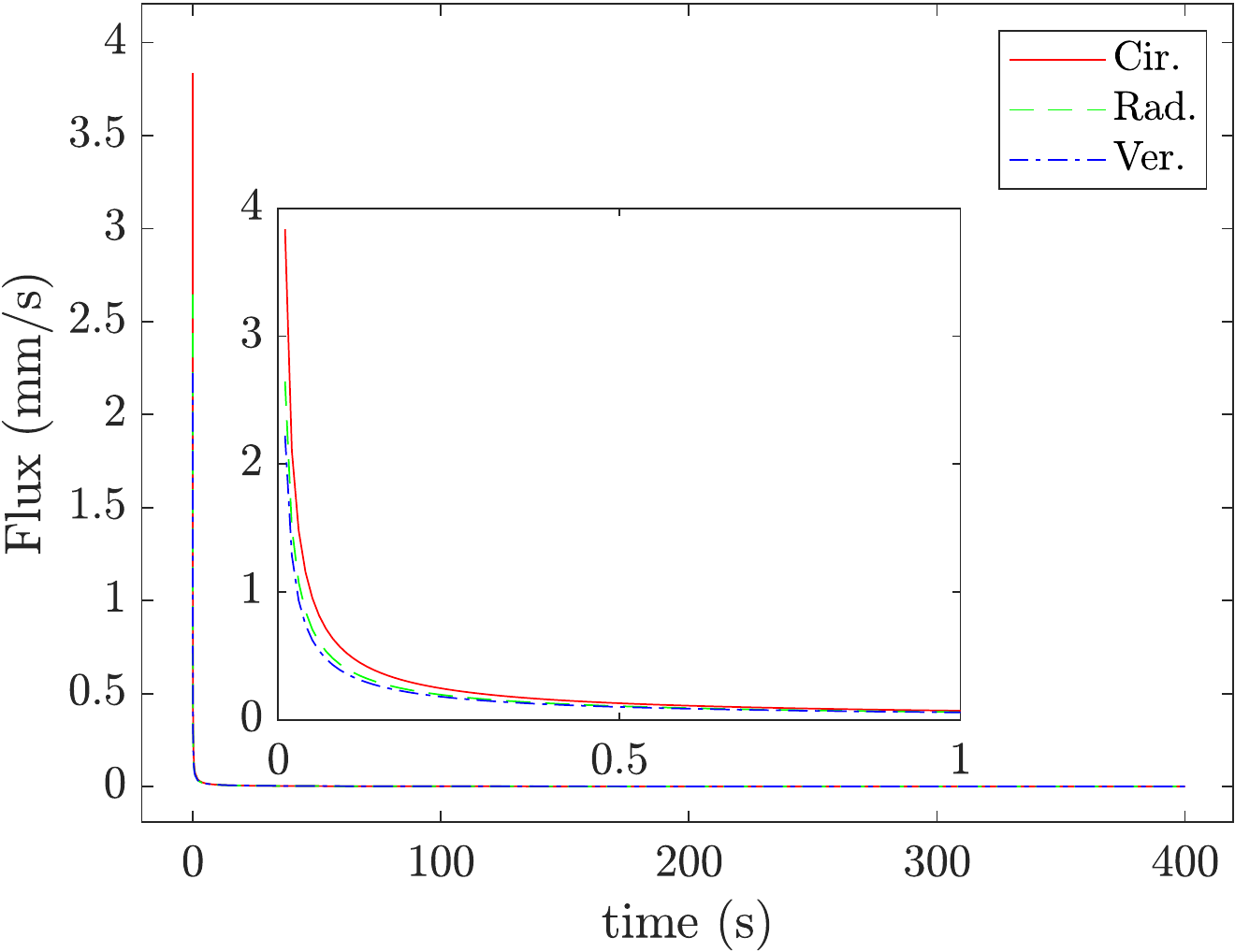}
\caption{}
\label{sf:Flux_directional_compar}
\end{subfigure}%
\caption{Flux, displacement and pore pressure for different directions (Circumferential, Radial and Vertical)}
\label{f:DirectionalCompar}
\end{figure}

To study the anisotropic behaviour of the meniscus, the confined creep test is performed as described in \sref{sec:ConfinedCompression} with material properties taken from the average values of the body region in circumferential, radial and vertical directions. For the numerical study, the height of the sample is taken as 3mm and pressure is increased to 0.07 MPa from zero and held constant. \fref{f:DirectionalCompar} shows the pore pressure, displacement and the flux out of the sample from the bottom (computed using \eref{eq:fluxfullform}) as a function of time for the three different regions.

\subsection{Numerical modelling of fractional consolidation} \label{AbaqusCreep}
\begin{figure}[htpb]
\centering
\includegraphics[width=0.3\textwidth]{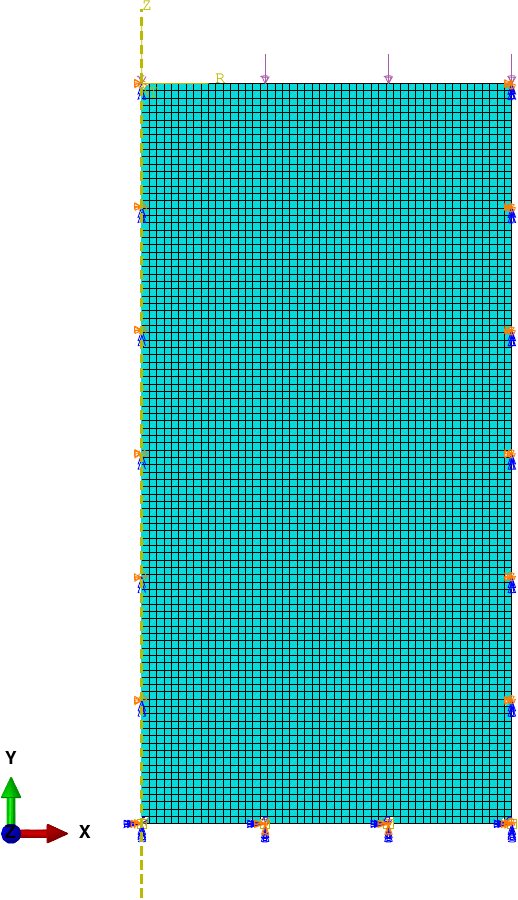}
\caption{Axisymmetric model employed for the confined creep numerical simulation}
\label{fig:AxisSymm_BCs}
\end{figure}

Numerical modelling of confined creep experiments of the meniscus is done in Abaqus using the parameters obtained from fitting the data given in \tref{T:MaterialProperties}. Poisson's ratio and Young's modulus are the elastic parameters required for Abaqus. As the confined compression tests depend only on the aggregate modulus, Poisson's ratio is arbitrarily assumed to be 0.3, and Young's modulus is calculated using the aggregate modulus and Poisson's ratio. Owing to symmetry, an axisymmetric model is considered for the numerical simulation. \fref{fig:AxisSymm_BCs} shows an axisymmetric model employed for this study. The domain is discretized with 4-noded bilinear quadrilateral elements (CAX4T) and a structured mesh is used. Based on a systematic mesh convergence study, a mesh size of 0.05mm was found to be adequate to model the behaviour. Axis-symmetric boundary conditions are enforced on the left of the computational domain. The displacements are restrained on the right and bottom faces. Zero pore pressure condition is applied at the bottom for free fluid flow and pressure $P_A=0.07MPa$ is applied at the top surface as a step load. For the numerical simulation, a fixed time increment of $\Delta t=$ 0.1~s is used, and the simulation is carried over for 400~s.

A comparison between the numerical, theoretical and experimental results is shown in \fref{fig:CreepDispCompar} and it is inferred that a good agreement is seen between the numerical, theoretical and experimental results. For the subsequent studies, the numerical models are used for complex loading in further sections.

\subsection{Fractional poroelastic model validation}
To further validate the proposed poroelastic model, a confined compression creep test with and without initial ramp loading and confined compression stress relaxation is done numercally and compared with experiments. We also compute the weight loss from the confined compression creep test. 

\begin{figure}[H]
            \centering
            \begin{subfigure}{.5\textwidth}
                \centering
                \includegraphics[width=\linewidth]{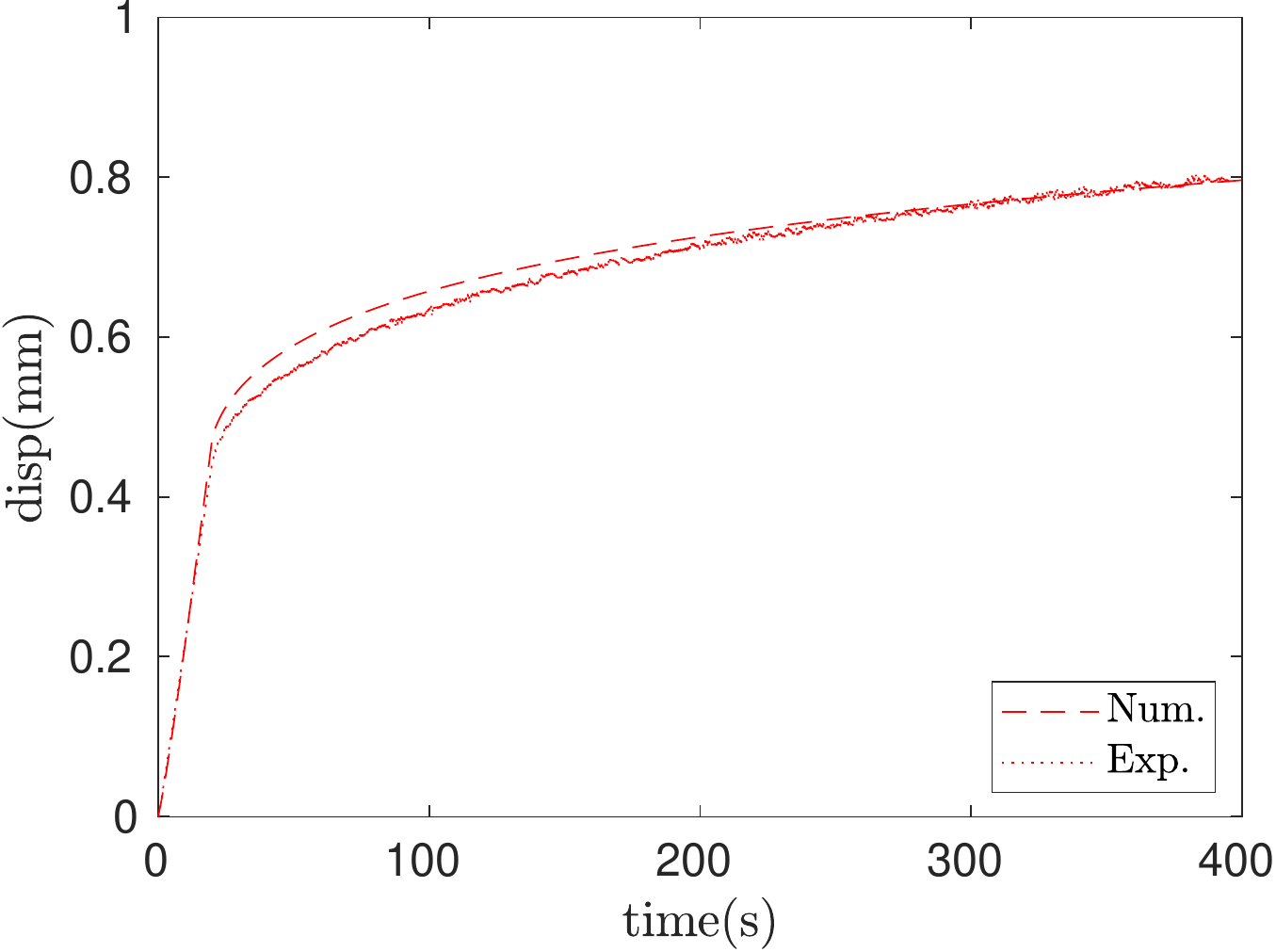}
                \caption{TK11 Body Circumferential}
                \label{TK11BC_CwR_AbqExp}
            \end{subfigure}%
            \begin{subfigure}{.5\textwidth}
                \centering
                \includegraphics[width=\linewidth]{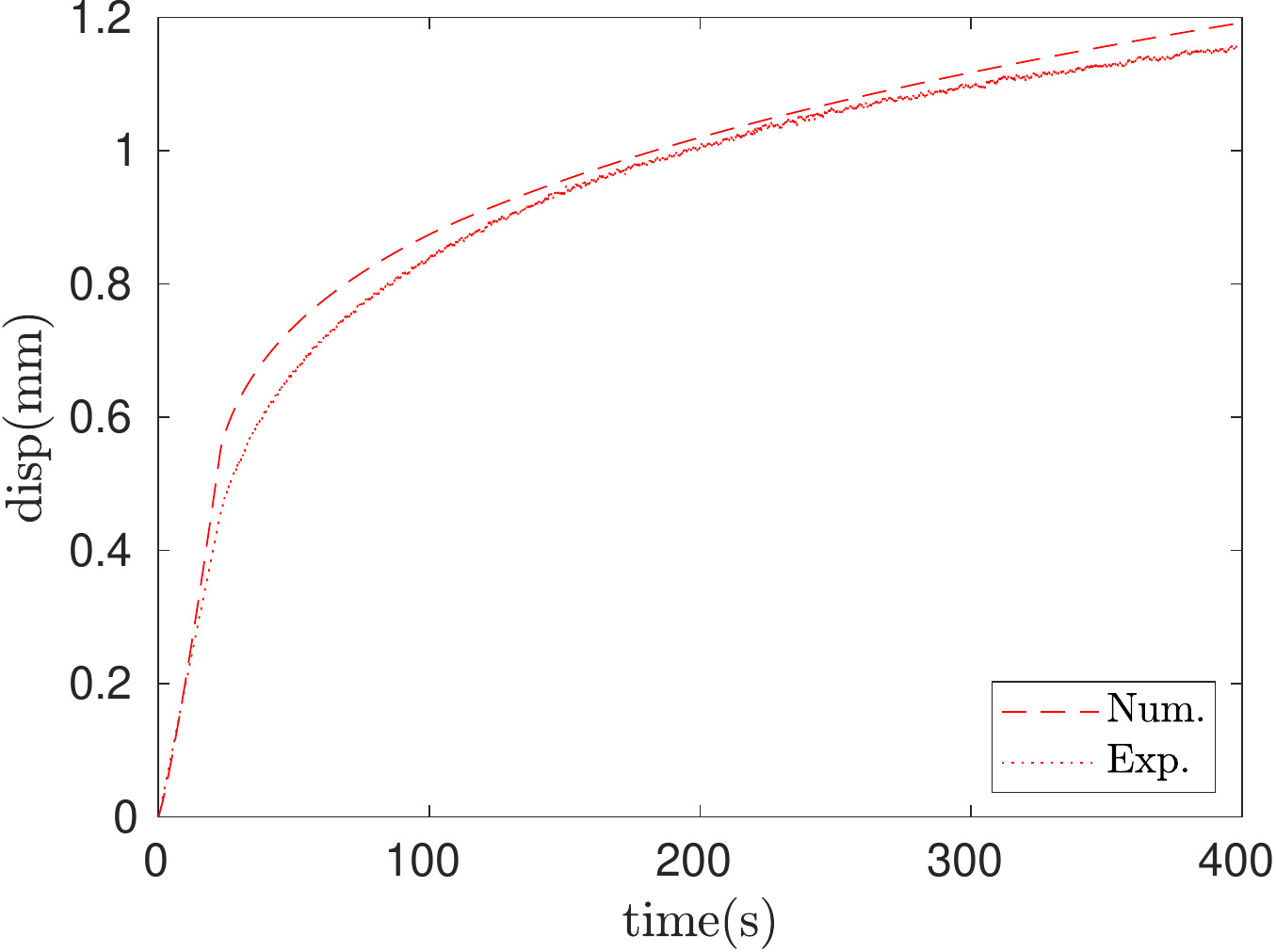}
                \caption{TK11 Body Radial}
                \label{TK11BR_CwR_AbqExp}
            \end{subfigure}\\%
            \begin{subfigure}{.5\textwidth}
                \centering
                \includegraphics[width=\linewidth]{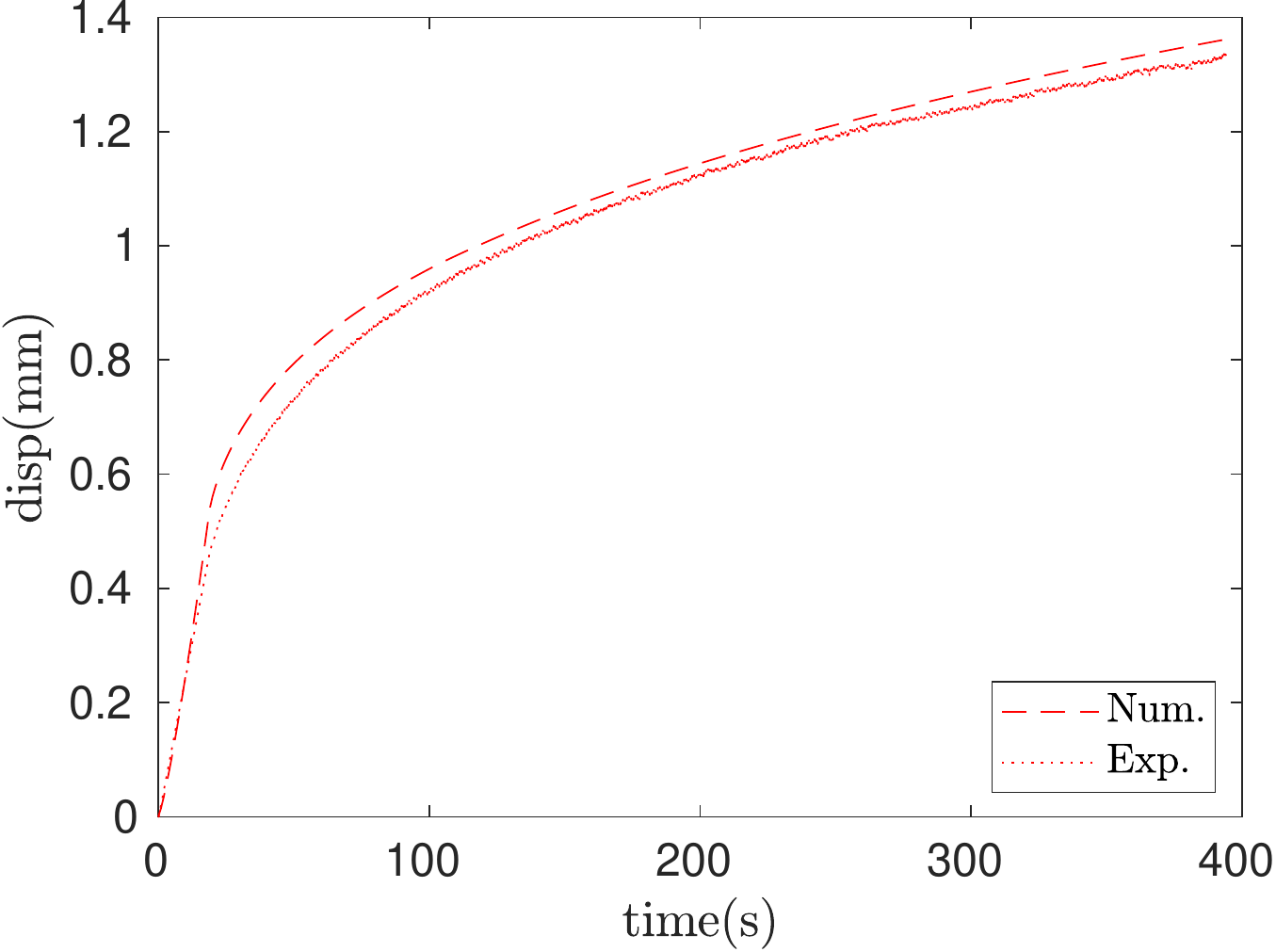}
                \caption{TK11 Body Vertical}
                \label{TK11BV_CwR_AbqExp}
            \end{subfigure}%
            \begin{subfigure}{.5\textwidth}
                \centering
                \includegraphics[width=\linewidth]{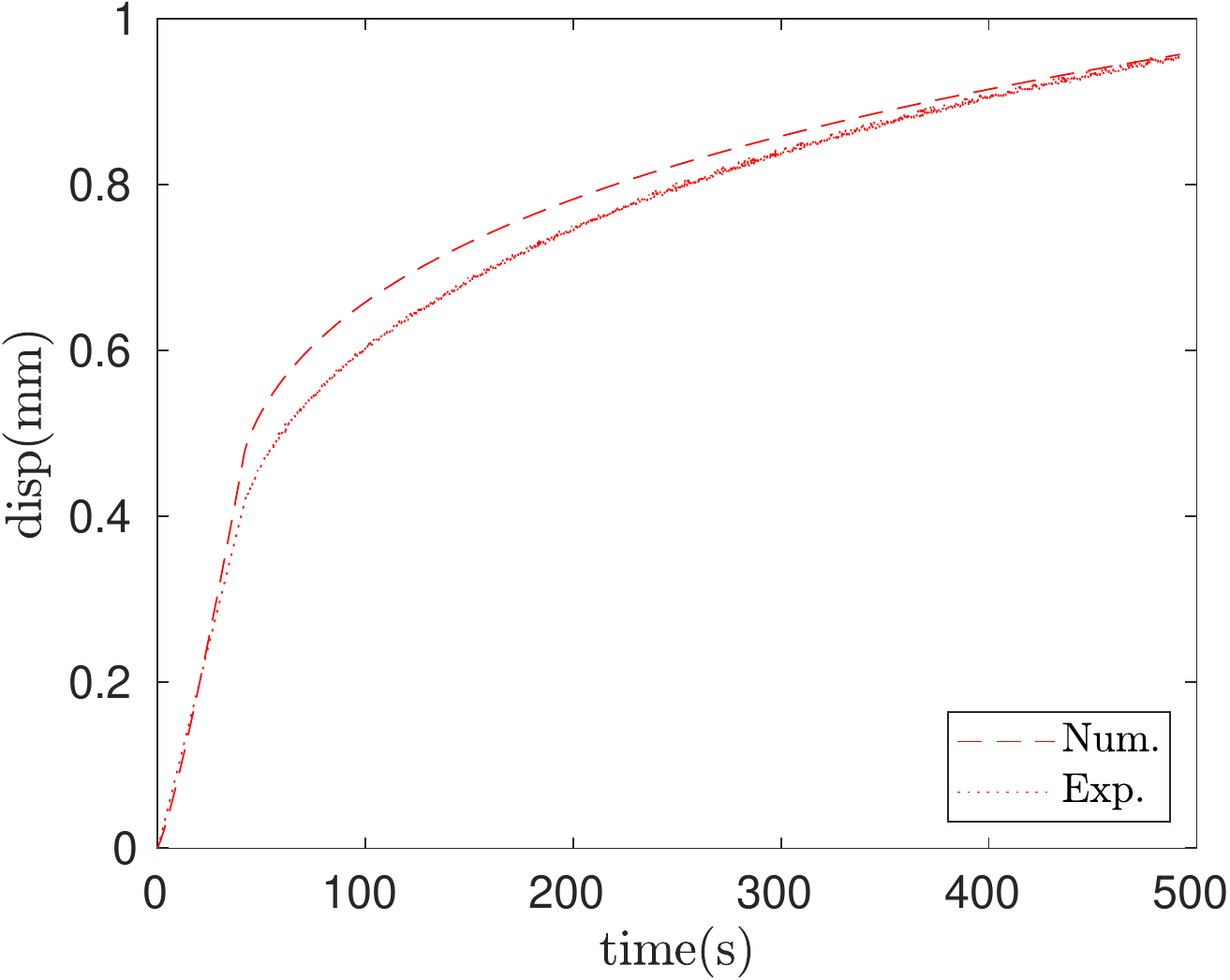}
                \caption{TK16 Anteriorhorn Radial}
                \label{TK16AR_CwR_AbqExp}
            \end{subfigure}%
            \caption{Validation by a creep with ramp test (a) Sample taken from the body region of TK11 sample in the circumferential direction, (b) Sample taken from the body region of TK11 sample in the radial direction, (c) Sample taken from the body region of TK11 sample in the vertical direction, (d)Sample taken from the anterior horn region of TK16 sample in the radial direction.}
            \label{fig:CreepWR}
        \end{figure}

\paragraph{Confined compression creep test with initial ramp loading}  
\label{sec:CwRCompar}
In the earlier study, the initial ramp present in the experiments was assumed as a step load. However, to assess the impact of this assumption on the model fitting parameters, we ran the FE model with the initial ramp (similar to that in the experiments) for the creep test and the results from the numerical simulation are compared with experiments. The creep test with ramp loading within Abaqus is modelled by providing the amplitude of the load in the load module with a ramp for a specific time period similar to experiments and then kept constant for the rest of the simulation. \fref{fig:CreepWR} shows the comparison of numerical results with the experiments for a few samples, in particular for samples TK11BC, TK11BR, TK11BV and TK16AR. The results of creep with ramp in \fref{fig:CreepWR}
show slight deviation at the end of ramp loading caused due to this assumption. But the deviation is small and can be neglected.
Therefore we conclude that the material properties extraction process as explained above is correct.

\begin{figure}[H]
        \centering
        \begin{subfigure}{.5\textwidth}
            \centering
            \includegraphics[width=\linewidth]{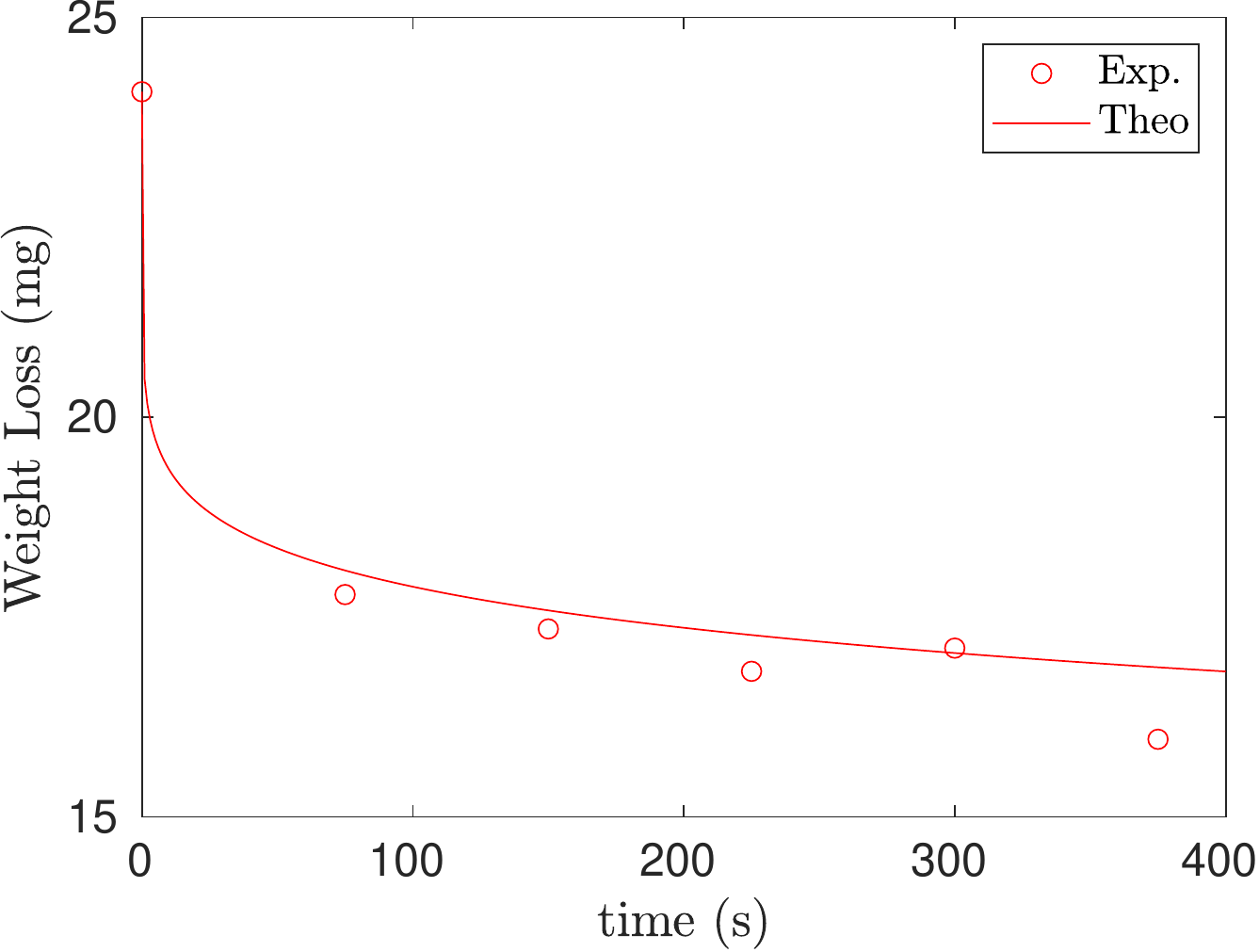}
            \caption{TK16 Body Circumfrential}
            \label{TK16BC_Wt_ExpTh}
        \end{subfigure}%
        \begin{subfigure}{.5\textwidth}
            \centering
            \includegraphics[width=\linewidth]{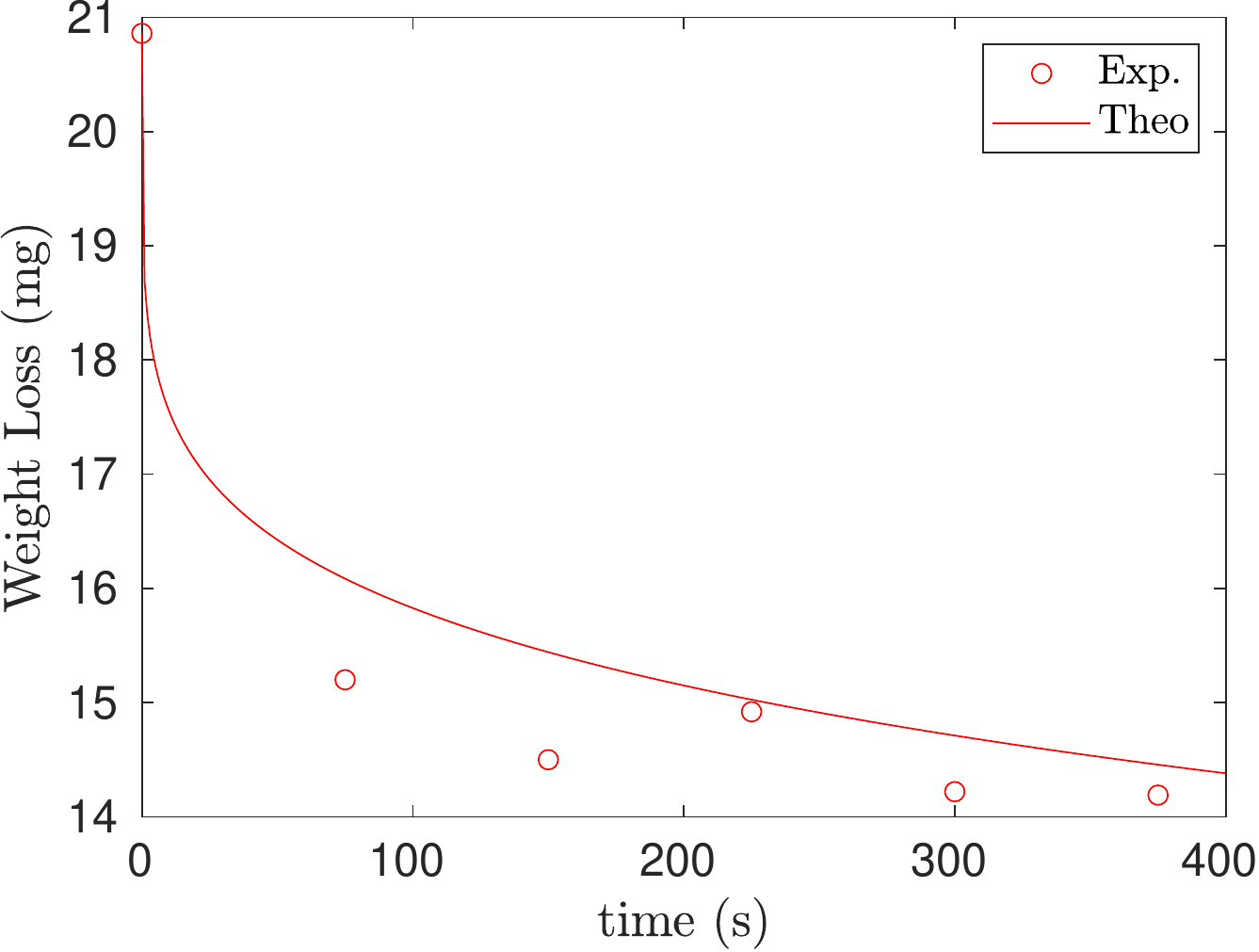}
            \caption{TK16 Body Radial}
            \label{TK16BR_Wt_ExpTh}
        \end{subfigure}\\%
        \begin{subfigure}{.5\textwidth}
            \centering
            \includegraphics[width=\linewidth]{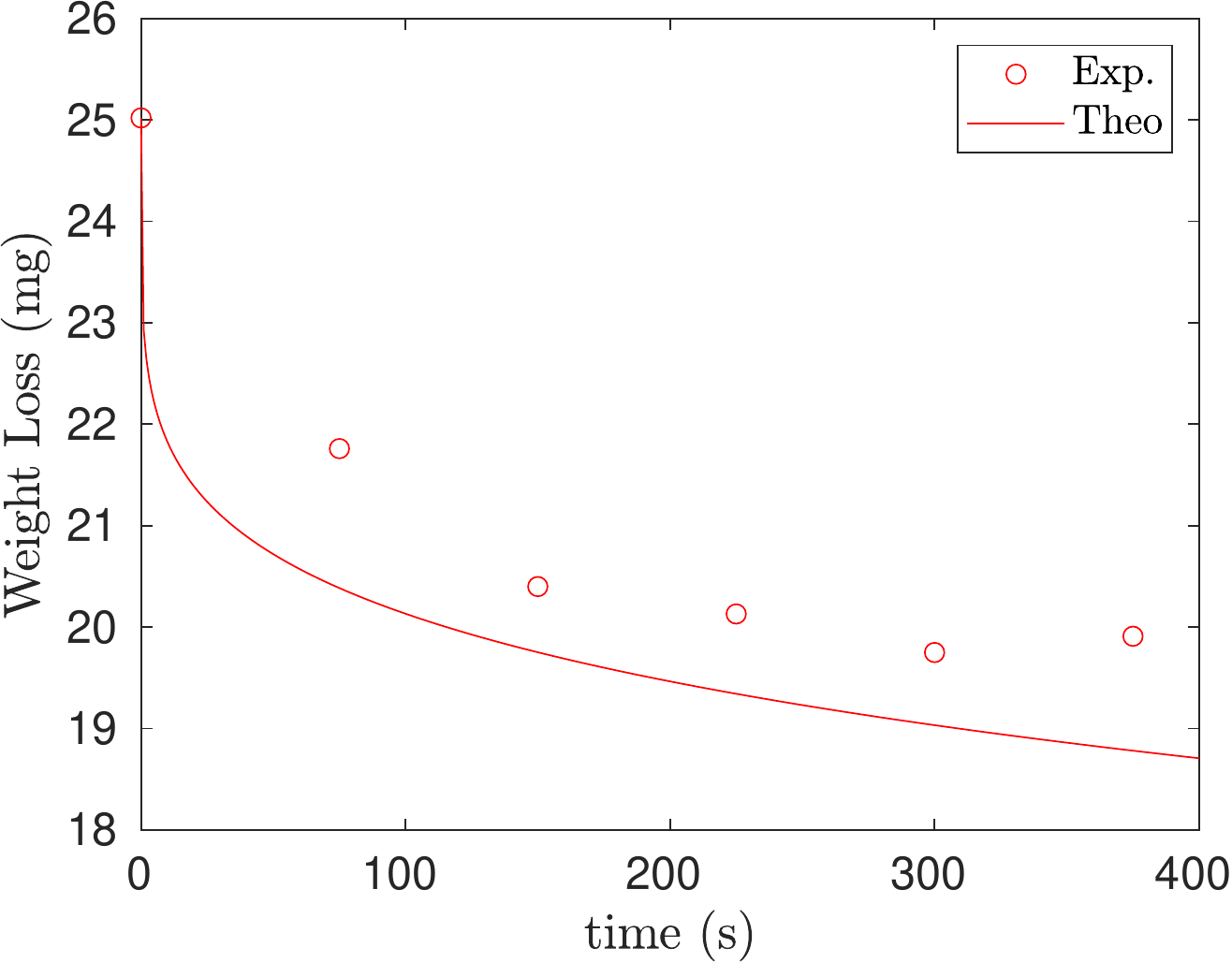}
            \caption{TK16 Body Vertical}
            \label{TK36BC_Wt_ExpTh}
        \end{subfigure}%
        \begin{subfigure}{.5\textwidth}
            \centering
            \includegraphics[width=\linewidth]{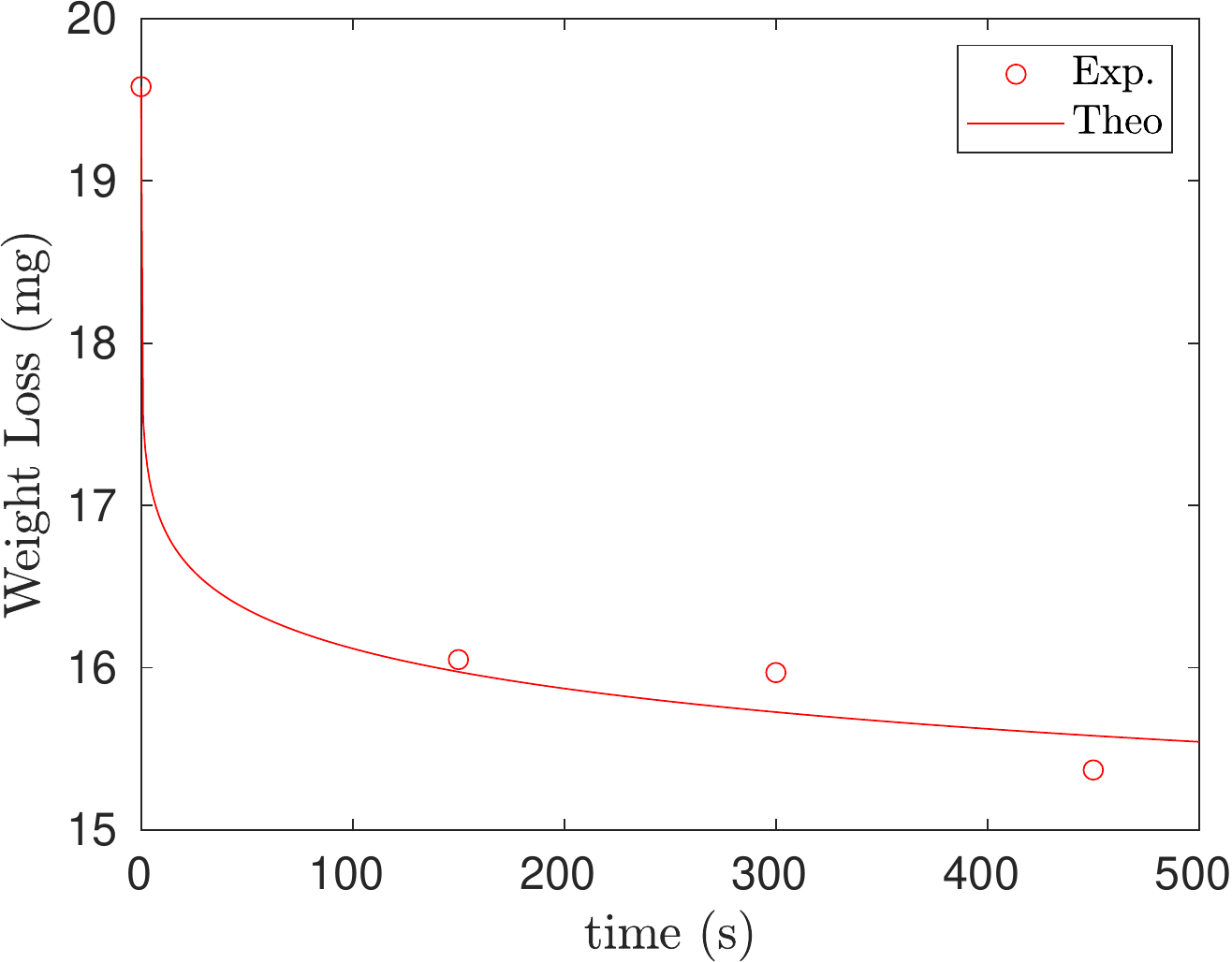}
            \caption{TK17 Anteriorhorn Radial}
            \label{TK17AR_Wt_ExpTh}
        \end{subfigure}\\%
        \caption{Comparison of theoretical Weight loss from the sample using parameters obtained from fitting with the experiments (a) Sample taken from the body region of TK16 sample in the circumferential  direction, (b) Sample taken from the body region of TK16 sample in the radial direction, (c) Sample taken from the body region of TK16 sample in the vertical direction, (d) Sample taken from the Anterior horn region of TK17 sample in the radial direction.}
        \label{WtLossCompar}
    \end{figure}

\begin{figure}[H]
\centering
\begin{subfigure}{.5\textwidth}
\centering
\includegraphics[width=\linewidth]{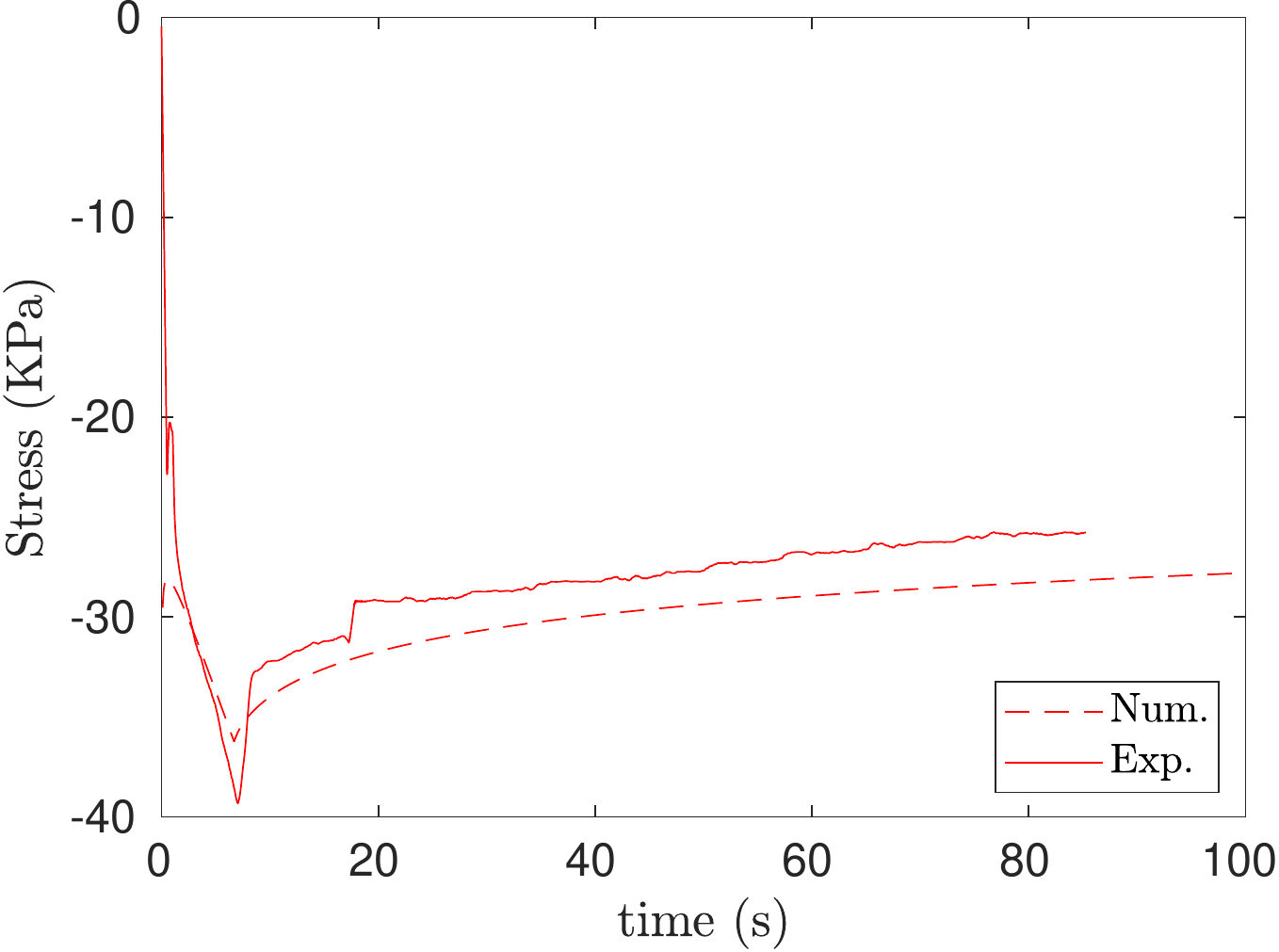}
\caption{TK11 Body Circumferential}
\label{TK11BC_SR_AbqExp}
\end{subfigure}%
\begin{subfigure}{.5\textwidth}
\centering
\includegraphics[width=\linewidth]{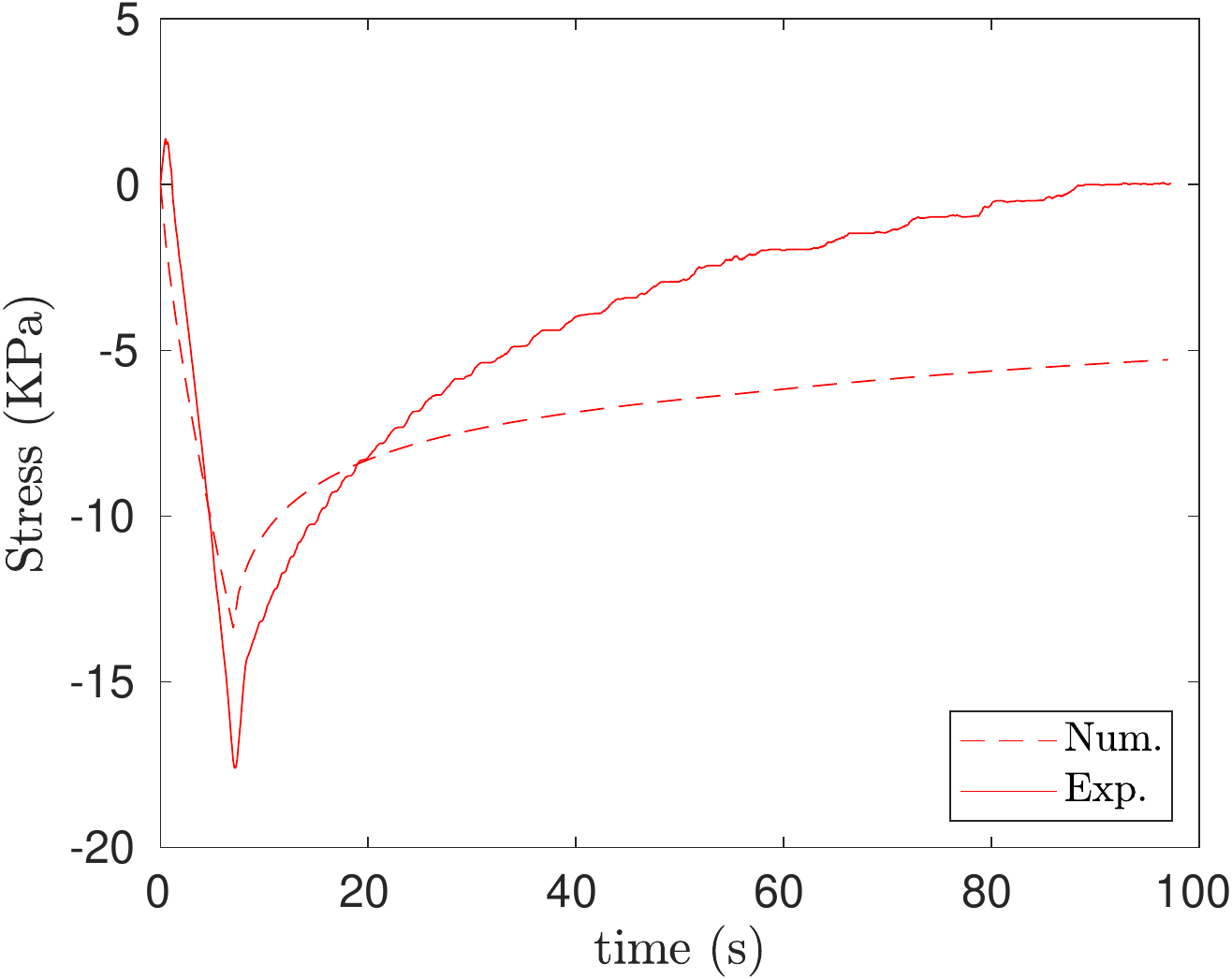}
\caption{TK11 Body Vertical}
\label{TK11BV_SR_AbqExp}
\end{subfigure}\\%
\begin{subfigure}{.5\textwidth}
\centering
\includegraphics[width=\linewidth]{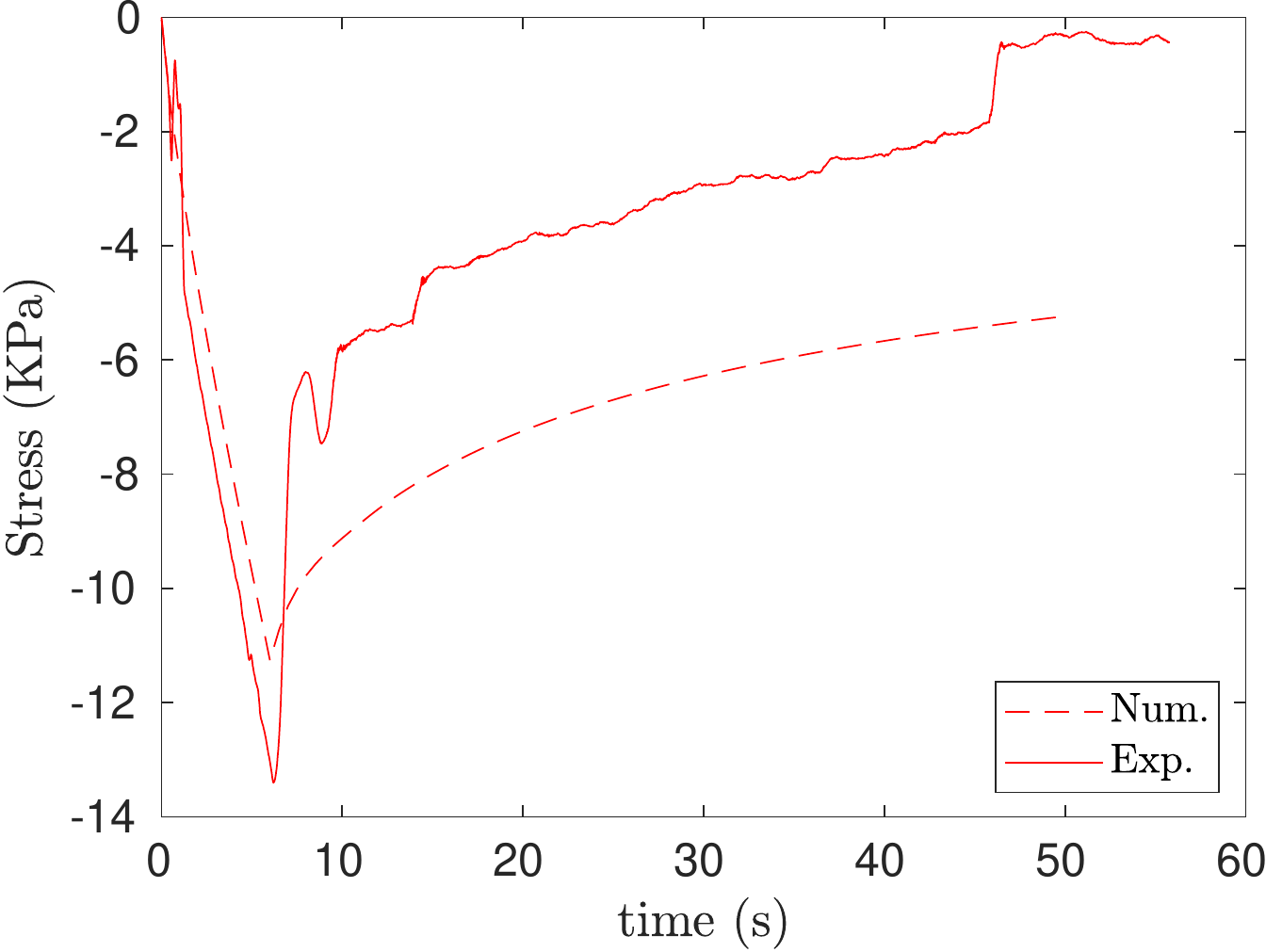}
\caption{TK11 Body Radial}
\label{TK11BR_SR_AbqExp}
\end{subfigure}%
\caption{Validation by comparing stresses from Abaqus with stress relaxation test from experiments (a) Sample taken from the body region of TK11 sample in the circumferential direction, (b) Sample taken from the body region of TK11 sample in the vertical direction, (c) Sample taken from the body region of TK11 sample in the radial direction.}
\label{fig:StrssRlxn}
\end{figure}

\paragraph{Weight loss from confined compression creep test}
\label{sec:WtLossCompar}
Weight loss for confined creep compression problem presented in  \sref{sec:ConfinedCompression} and compared with the experiments carried out by Bulle et al.\citep{Bulle2021}. Specific weight of water at 37$^\circ$C is used($w_s =$ 997Kg/m$^3$). Fluid flow is collected at the bottom of the cylindrical sample with the area of the circular cross-section $\boldsymbol{A}$ with diameter $d=3mm$. Material properties are taken from the \tref{T:MaterialProperties} which were obtained from fitting the fractional poroelastic model with displacement data from the confined creep test. Weight loss with time for some samples is shown in the \fref{WtLossCompar}.
Theoretical and experimental results for a few samples match as shown in the \fref{WtLossCompar}. Moreover, they do not match well for some samples, as shown in \fref{WtLossCompar2}. This could be due to the measurement error as it was carried out manually. The weight of the sample is measured after every 75 seconds of creep by taking it out.
\paragraph{Confined compression stress relaxation}
\label{sec:SRCompar}
The fractional poroelastic model is validated using the confined compression stress relaxation tests performed by Bulle et al.~\citep{Bulle2021}. Tests were performed on the cylindrical sample of height $h$, confined on the sides, and the fluid can flow freely from the bottom. Tests consist of preconditioning the sample with 10\%$h$ compression at the ramp velocity 0.3\%$h$/s. Then, five relaxation steps with 2\%$h$ compression at a ramp velocity of 0.3\%$h$/s s were performed. In all the cases, the force applied is measured with time. The finite element modelling of this confined compression stress relaxation test is performed in Abaqus similar to the creep test as mentioned in  \sref{AbaqusCreep}. Material properties obtained from fitting the creep data were used. In this test, the displacement is given as input, and the force is measured as output. \fref{fig:StrssRlxn} shows a comparison of the FEM results 
with the experiments. The comparisons of the model with the experiments are reasonably good. 

\section{Discussion}
We show that the fractional poroelastic model is more appropriate to model the time-dependent behaviour of soft tissue presenting a hierarchically arranged porous architecture. The model has been verified and validated. Our parameter fittings show that the RMS error = 1.42$\times$10$^{-5}$ for the fractional poroelastic model is considerably lower than the one given by the classical model that had an RMS error = 6.53$\times$10$^{-5}$, see \fref{fig:ClassicvsFrac}. The theoretical displacement curve starts from 0, whereas the experiment shows an initial displacement. This is attributed to the assumption of undrained conditions at the initial time and the assumption of incompressible material. We obtain a value of fractional order $\beta$ in the range of $\beta=$ 0.51--0.73 as shown in \tref{Meniscus param sample}. From the numerical study, it can be inferred that the transport phenomenon during the confined compression is ruled by a fractional version of Darcy's law. 
At the beginning of the confined compression test, the solid is fully saturated, therefore the pore pressure reaches a maximum value (see \fref{fig:pvst_AbqTh_comp}). This is also predicted by Terzaghi's consolidation theory, in line with the fact that the fluid pressure entirely carries the load. During the consolidation process, the fluid flows out of the sample at a rate depending on the anomalous permeability and on the order of the fractional derivative (see \fref{fig:AbqTh_comp}). The higher the value of $\beta$ the faster the diffusion process. As the fluid flows out of the sample, the pore pressure decreases and the solid starts deforming (see \fref{fig:dvst_AbqTh_comp}). Also, higher values of $\beta$ imply a faster solid deformation process.
We notice that the aggregate modulus $M$ obtained from fitting mentioned in \tref{T:MaterialProperties} is of the same order as the literature \citep{Sweigart2004, John2013, Andreas2013, Helena2008}. The anomalous permeability $\lambda_{\beta}$ is in the order of $10^{-12}/10^{-13}$. From the statistical analysis, it can be concluded that the meniscus in the body region is elastically isotropic, and pore pressure diffusion is transversely isotropic with symmetry in the vertical and radial directions. Higher hydraulic permeability in the circular direction can be attributed to the fact that the fibres are oriented in a circular direction in the body region, see \fref{fig:RegionsnDirections}. To model the anisotropy in the anomalous transport phenomena, the following anisotropic form of fractional Darcy's law is proposed:
\begin{equation}
\mathbf{J}_p =
\underbrace{\begin{bmatrix}
\lambda_{rad} & 0 & 0\\
0 & \lambda_{cir} & 0\\
0 & 0 & \lambda_{ver}
\end{bmatrix}}_\mathbf{\Lambda}
\underbrace{\begin{bmatrix}
\prescript{}{0}D_t^{\beta_{rad}} & 0 & 0\\
0 & \prescript{}{0}D_t^{\beta_{cir}} & 0\\
0 & 0 & \prescript{}{0}D_t^{\beta_{ver}}\\
\end{bmatrix}}_{\mathbf{\prescript{}{0}D_t^\mathbf{\beta}}}
\begin{bmatrix}
\nabla p
\end{bmatrix}
\label{eq:AnisoFluidFlux}
\end{equation}
where, $\mathbf{\Lambda}$ is the permeability tensor and $\mathbf{\prescript{}{0}D_t^\mathbf{\beta}}$ is the fractional derivative operator. We notice that for the body region,  the anomalous permeability tensor $\mathbf{\Lambda}$ is transversely isotropic with $\lambda_{\rm cir}$ considerably higher than $\lambda_{\rm ver}$ and $\lambda_{\rm rad}$. Similarly $\mathbf{\beta}$ tensor is transversely isotropic with $\beta_{\rm cir}$ considerably higher than the $\beta_{\rm rad}$ and $\beta_{\rm ver}$. Consequently, from the \fref{sf:Flux_directional_compar}, it is observed that the fluid flow rate through the sample's base during the consolidation experiment  in the circumferential direction is higher than the other two directions.
Pore pressure in the circumferential direction reaches a steady state faster than the other two directions, as shown in the \fref{sf:porep_directional_compar} because of the higher permeability in the circumferential direction.
Based on the mean values and the standard deviation in the anterior horn region, it is opined that the aggregate modulus in the radial direction is higher than the other two directions, viz., vertical and circular directions. It is inferred that the permeability and the order of derivative are higher in the radial and vertical directions than in the circumferential directions. The elastic tensor results are transversely isotropic, with the modulus in the radial direction being higher, and both $\mathbf{\Lambda}$  and $\mathbf{\beta}$ are transversely isotropic with circumferential being lower. Different properties in different directions show that the elastic, anomalous permeability and the order of the derivative tensors are transversely isotropic. Furthermore, properties vary with the anterior horn and body region, showing that it is not homogeneous. 
\begin{figure}[H]
    \centering
    \includegraphics[width=\textwidth]{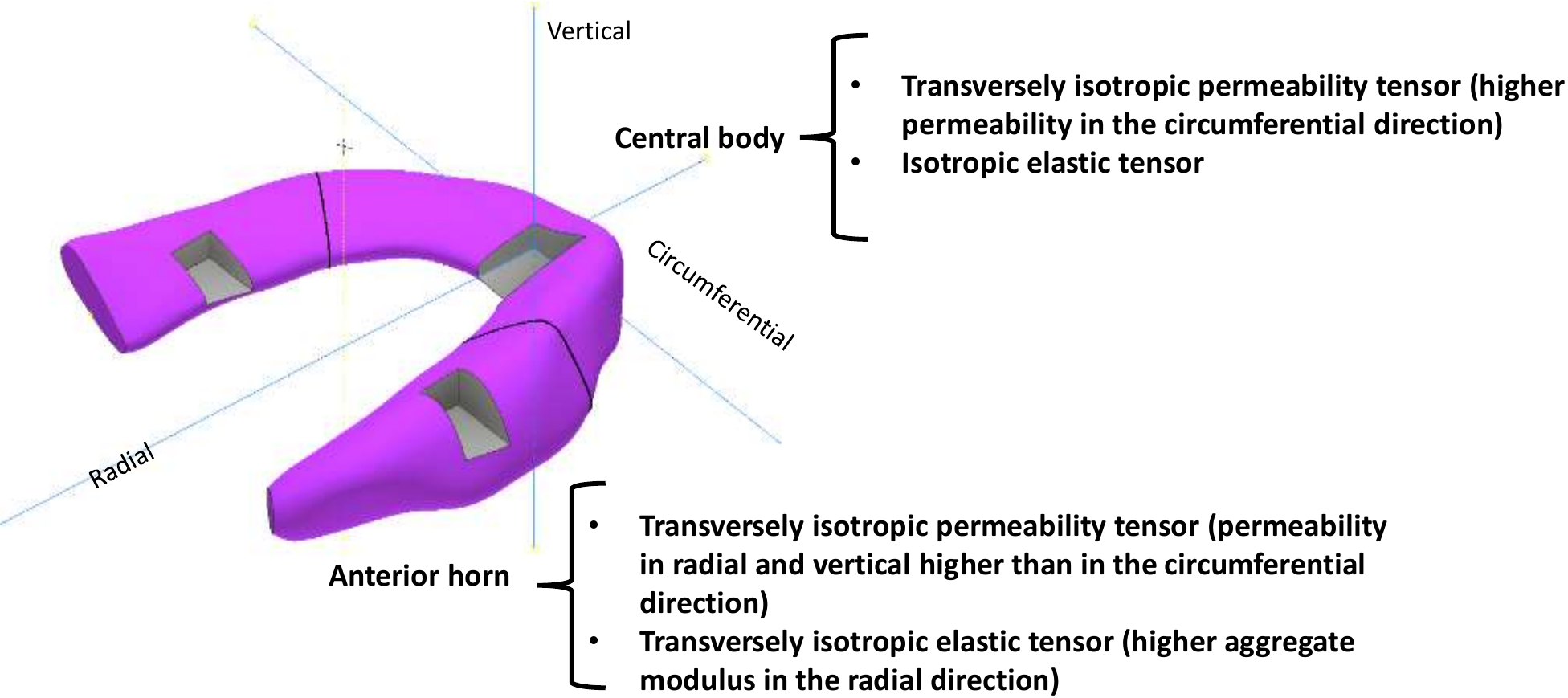}
    \caption{Directional and regional dependence of material properties in the meniscus.}
    \label{fig:meniscus_anisotropy}
\end{figure}
In summary (see \fref{fig:meniscus_anisotropy}), we show that the body region which has a load-bearing function exhibits a transversely isotropic behaviour due to the rate of fluid flow being about three times higher (faster diffusion) in the circumferential direction which is consistent with the preferential direction of collagen channels. It explains the role of fluid pressure in sustaining the load. Furthermore, we show that both the elastic and the permeability tensors are transversely isotropic in the anterior horn. The aggregate modulus is higher in the radial direction compared to the circumferential and vertical directions. It explains the role of the anterior horn to be compliant in the circumferential and radial directions to accommodate the kinematics of the tissue.

\section{Conclusions}
The meniscus has a porous, hierarchical, multi-oriented fiber and bundle structure filled with fluid that provides optimized load support and lubrication properties. To understand its function, confined compression creep tests were performed on different regions of meniscus tissue and in different orientations. Biot's theory with fractional Darcy's law 
is used to find the material properties. It is observed that, with the classical Darcy's law, fitting gives an RMS error of 6.53$\times$10$^{-5}$, while the fractional Darcy's law gives an RMS error of 1.42$\times$10$^{-5}$. It is shown that Biot's theory with fractional Darcy's law is better suited for modelling the poroelastic behaviour of the meniscus. Three material parameters, viz., aggregate modulus($M$), fractional order($\beta$) and the permeability($\lambda_\beta$) required for the fractional Biot's theory were found by fitting the theory with the confined compression creep experiment results for all the samples. Material properties in different directions of the body region are compared using ANOVA test methods. It is observed that the aggregate modulus ($M_{\rm cir} =$ 75.4 $ \pm$48.9 KPa, $M_{\rm rad} =$ 63.7$ \pm$14.1  KPa and $M_{\rm ver} =$ 65.3$ \pm$13.7 KPa) is isotropic and the permeability ($\lambda_{\beta_{\rm cir}} =$ 3.75$\times$10$^{-12}\pm$ 2.36$\times$10$^{-12}$m$^2$/Pa.s$^{1-\beta}$, $\lambda_{\beta_{\rm rad}} =$ 1.25$\times$10$^{-12}\pm$ 4.04$\times$10$^{-13}$m$^2$/Pa.s$^{1-\beta}$ and $\lambda_{\beta_{\rm ver}} =$ 8.61$\times$10$^{-13}\pm$ 1.94$\times$10$^{-13}$m$^2$/Pa.s$^{1-\beta}$) and the fractional order ($\beta_{\rm cir} =$ 0.73$\pm$0.00, $\beta_{\rm rad} = 0.64\pm 0.06$ and $\beta_{\rm ver} =$ 0.58$\pm$ 0.04) are transversely isotropic with properties in the circumferential direction is greater than the vertical and radial directions.
    
Biot's theory with fractional Darcy's law is implemented numerically using the FEM in Abaqus software using UMATHT and UMAT subroutines. This fractional poroelastic theory is validated by using the material properties obtained from fitting confined creep tests to model stress relaxation, creep with ramp, and weight loss tests. The results are then compared with experimental results. The comparison shows that Biot's theory with fractional Darcy's law is capable of better characterising the poroelastic behaviour of the meniscus. Flux out of the meniscus, which is difficult to measure experimentally, can easily be computed efficiently using the current model. In future works, Biot's theory with fractional Darcy's law will be extended to include anisotropy and use it to understand the working of the meniscus in the knee joint. 
    
\section{Acknowledgements}
O.B would like to acknowledge the European Union’s Horizon 2020 -EU.1.3.2. - Nurturing excellence by means of cross-border and cross-sector mobility under the Marie Skłodowska-Curie individual fellowship MSCA-IF-2017, MetaBioMec, Grant agreement ID: 796405. The authors thank the Rizzoli Orthopaedic Institute and, in particular, G. Marchiori and M. Berni for their invaluable insights on the experiments. 
\bibliographystyle{elsarticle-num}
\bibliography{references}






\appendix
\section{}
Using \eref{eq:Strnpp}, the displacement field for the problem described in \sref{sec:ConfinedCompression} can be solved as follows:
\begin{equation}
\epsilon_{zz} = \dfrac{\partial u_z}{\partial z} = \dfrac{3}{3K+4G}\Big(-P_A +\alpha p\Big) 
\label{eq:DispPorep}
\end{equation}
Substituting \eref{eq:PorepSoln} in \eref{eq:DispPorep} and integrating
\begin{equation}
u(z,t)=\dfrac{3P_{A}}{3K+4G}\Bigg[-z+\gamma \alpha \sum \limits_{n=1,3}^\infty  E_{1-\beta,1} \left(-{\dfrac{n^{2}\pi^{2} \bar{\lambda} t^{1-\beta}}{4h^{2}}}\right) b_{n}\sin \dfrac{n \pi z}{2h} \Bigg] + d
\label{eq:10a}
\end{equation}
where, $b_{n}= \frac{8h}{(n \pi)^{2}} \left(-1\right)^{\frac{n-1}{2}}$. Applying the boundary condition given by \sref{sec:governing eqns}, we get:
\begin{equation}
u(h,t) = \frac{3P_{A}}{3K+4G}\Bigg[-h+\gamma \alpha \sum \limits_{n=1,3}^\infty  E_{1-\beta,1} \left(-{\frac{n^{2}\pi^{2} \bar{\lambda} t^{1-\beta}}{4h^{2}}}\right) \frac{8h}{(n \pi)^{2}} \left(-1\right)^{\frac{n-1}{2}} \sin \frac{n \pi }{2} \Bigg] + d = 0
\end{equation}
Since, 
\begin{equation}
\Big(-1\Big)^{\frac{n-1}{2}} \sin \frac{n \pi }{2} = 1 \hspace{20pt} \text{for} \qquad  n = 1,3,\cdots 
\end{equation}  
Therefore, the constant $d$ is
\begin{equation}
d = -\dfrac{3P_{A}}{3K+4G}\Bigg[-h+\gamma \alpha \sum \limits_{n=1,3}^\infty  E_{1-\beta,1} \left(-{\frac{n^{2}\pi^{2} \bar{\lambda} t^{1-\beta}}{4h^{2}}}\right) \dfrac{8h}{(n \pi)^{2}} \Bigg]
\label{const_d}
\end{equation}
Substituting \eref{const_d} in \eref{eq:10a}, Displacement field equation is obtained
\begin{equation}
u(z,t)=\frac{3P_{A}}{3K+4G}\Bigg[(h-z)+\gamma \alpha \sum \limits_{n=1,3}^\infty  E_{1-\beta,1} \left(-{\frac{n^{2}\pi^{2} \bar{\lambda} t^{1-\beta}}{4h^{2}}}\right) \dfrac{8h}{(n \pi)^2}\bigg((-1)^{\frac{n-1}{2}}\sin \dfrac{n \pi z}{2h}-1 \bigg) \Bigg]
\label{eq:10a2}
\end{equation}

\section{}\label{AppendixB}
Weight loss from the sample over time W(t) for the problem described in \sref{sec:ConfinedCompression} can be solved using the fluid flux as:
\begin{equation}
W(t) = W_0 - \int_0^t\mathbf{J}_p\cdot w_s\mathbf{A}dt
\label{WeightLoss}
\end{equation}
where $\mathbf{J}_p$ is the fluid flux given by \eref{eq:FluidFlux}, $w_s$ is the specific weight, $\mathbf{A}$ is the cross-sectional area over which the flux is calculated and $W_0$ being the initial weight of the sample.
Differentiating \eref{eq:PorepSoln} with respect to $z$, we get:
\begin{equation}
\frac{\partial p }{\partial z } = -P_A \gamma \sum_{n=1,3}^{\infty} E_{1-\beta, 1}\bigg[-\frac{n^2 \pi^2 \bar{\lambda}t^{1-\beta}}{4h^2}\bigg]\dfrac{2}{h}(-1)^{\frac{n-1}{2}} \sin\bigg(\frac{n\pi z}{2h}\bigg)
\label{eq:PorepDeriv}
\end{equation}
Fluid can flow only at the bottom, as shown in \fref{sf:Schematic_Exp}. Hence, the weight loss is solved at $z=h$. \eref{eq:PorepDeriv} at $z=h$ is:
\begin{equation}
\dfrac{\partial p }{\partial z } = -P_A \gamma \sum_{n=1,3}^{\infty} E_{1-\beta, 1}\bigg[-\frac{n^2 \pi^2 \bar{\lambda}t^{1-\beta}}{4h^2}\bigg]\dfrac{2}{h}
\label{dpdzath}
\end{equation}
Using \eref{dpdzath} and \eref{eq:FluidFlux} in weight loss \eref{WeightLoss}, we get:
\begin{equation}
W(t) = W_0 - \lambda_\beta w_s A \prescript{}{0}D_t^{\beta - 1}\bigg[P_A \gamma \sum_{n=1,3}^{\infty} E_{1-\beta, 1}\bigg(-\frac{n^2 \pi^2 \bar{\lambda}t^{1-\beta}}{4h^2}\bigg)\dfrac{2}{h}\bigg]
\label{Wtloss2}
\end{equation}
Using the identity of the fractional derivative from \cite{Podlubny1999}, which is:
\begin{equation}
\prescript{}{0}D_t^\alpha t^{\beta - 1}E_{\mu, \beta}(\lambda t^\mu) = t^{\beta - 1 - \alpha} E_{\mu, \beta -\alpha}(\lambda t^\mu)
\label{FractDerIden}
\end{equation}
\eref{Wtloss2} can be found as:
\begin{equation}
W(t) = W_0 - \lambda_\beta w_s A P_A \gamma t^{1- \beta} \sum_{n=1,3}^{\infty} E_{1-\beta, 2-\beta}\bigg(-\dfrac{n^2 \pi^2 \bar{\lambda}t^{1-\beta}}{4h^2}\bigg)\dfrac{2}{h}
\end{equation}


\counterwithin*{table}{section}
\section{}
\label{A:MatProps_classic}
\begin{table}[H]
    \centering
\begin{tabular}{clccc}\toprule
S No. &Sample &M $\times 10^5(Pa)$ & $\lambda \times 10^{-13} (m^2/Pa.s) $ & RMS $\times 10^{-5}$ \\\midrule
1 &TK11BC &3.56 &6.95 &6.53 \\
2 &TK11BR &2.13 &3.91 &7.40 \\
3 &TK11BV &2.30 &4.44 &8.10 \\
4 &TK16BC1 &2.48 &5.32 &7.04 \\
5 &TK16BR2 &2.55 &2.45 &5.53 \\
6 &TK16BV &2.50 &3.00 &6.37 \\
7 &TK16AC1 &2.55 &2.80 &7.64 \\
8 &TK16AR1 &2.64 &2.42 &7.16 \\
9 &TK16AV &2.69 &1.79 &5.83 \\
10 &TK16BC2 &2.39 &8.87 &8.15 \\
11 &TK16BR1 &2.31 &2.68 &7.12 \\
12 &TK16BV2 &2.56 &2.74 &6.96 \\
13 &TK16PV &3.97 &12.1 &4.66 \\
14 &TK17BC &1.46 &11.8 &13.4 \\
15 &TK17BR &1.67 &4.17 &6.93 \\
16 &TK17BV &1.80 &1.99 &5.36 \\
17 &TK17AC &3.49 &208 &4.58 \\
18 &TK17AR &4.01 &23.3 &4.49 \\
19 &TK18BC &1.44 &1.33e4 &10.5 \\
20 &TK18BR &2.37 &5.27 &7.54 \\
21 &TK18BV &2.17 &2.89 &6.45 \\
22 &TK22BR &2.55 &4.66e3 &1.24 \\
23 &TK22AC &2.41 &3.69 &9.01 \\
24 &TK22AR &1.95 &4.74 &10.0 \\
25 &TK22AV &2.97 &1.31 &4.83 \\
26 &TK36BC &1.98 &9.33 &10.2 \\
27 &TK37BC &2.52 &16.5 &12.4 \\
28 &TK37BR &2.67 &3.23 &7.42 \\
29 &TK37BV &2.29 &2.66 &6.25 \\
\bottomrule
\end{tabular}
\label{t:MatProps_classic}
\caption{Material properties obtained from fitting using classical Biot's theory}
\end{table}

\counterwithin*{figure}{section}
\section{}
\begin{figure}[H]
\centering
\begin{subfigure}{.5\textwidth}
\centering
\includegraphics[width=\linewidth]{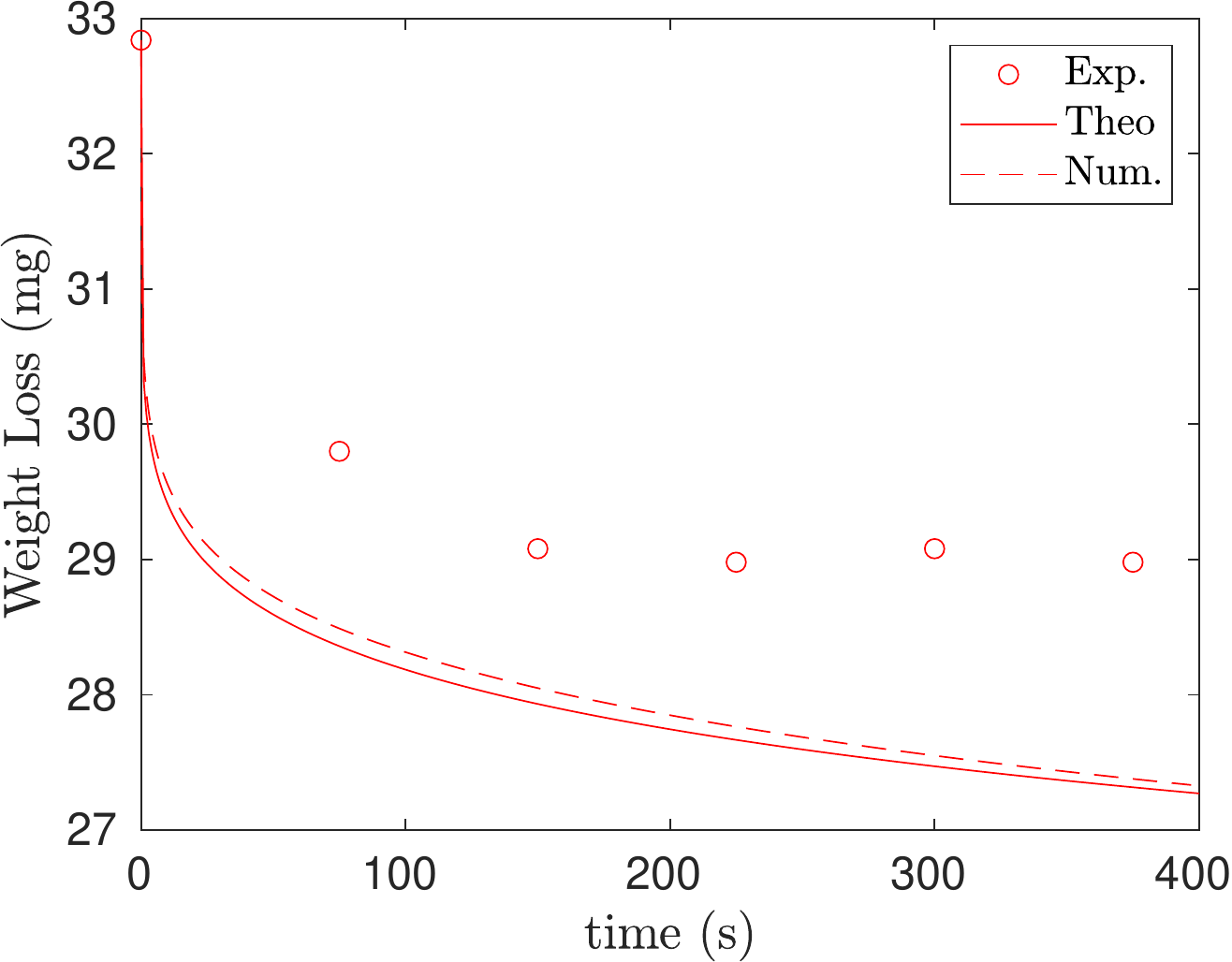}
\caption{TK11 Body Circ}
\label{TK11BC_Wt_AbqExpTh}
\end{subfigure}%
\begin{subfigure}{.5\textwidth}
\centering
\includegraphics[width=\linewidth]{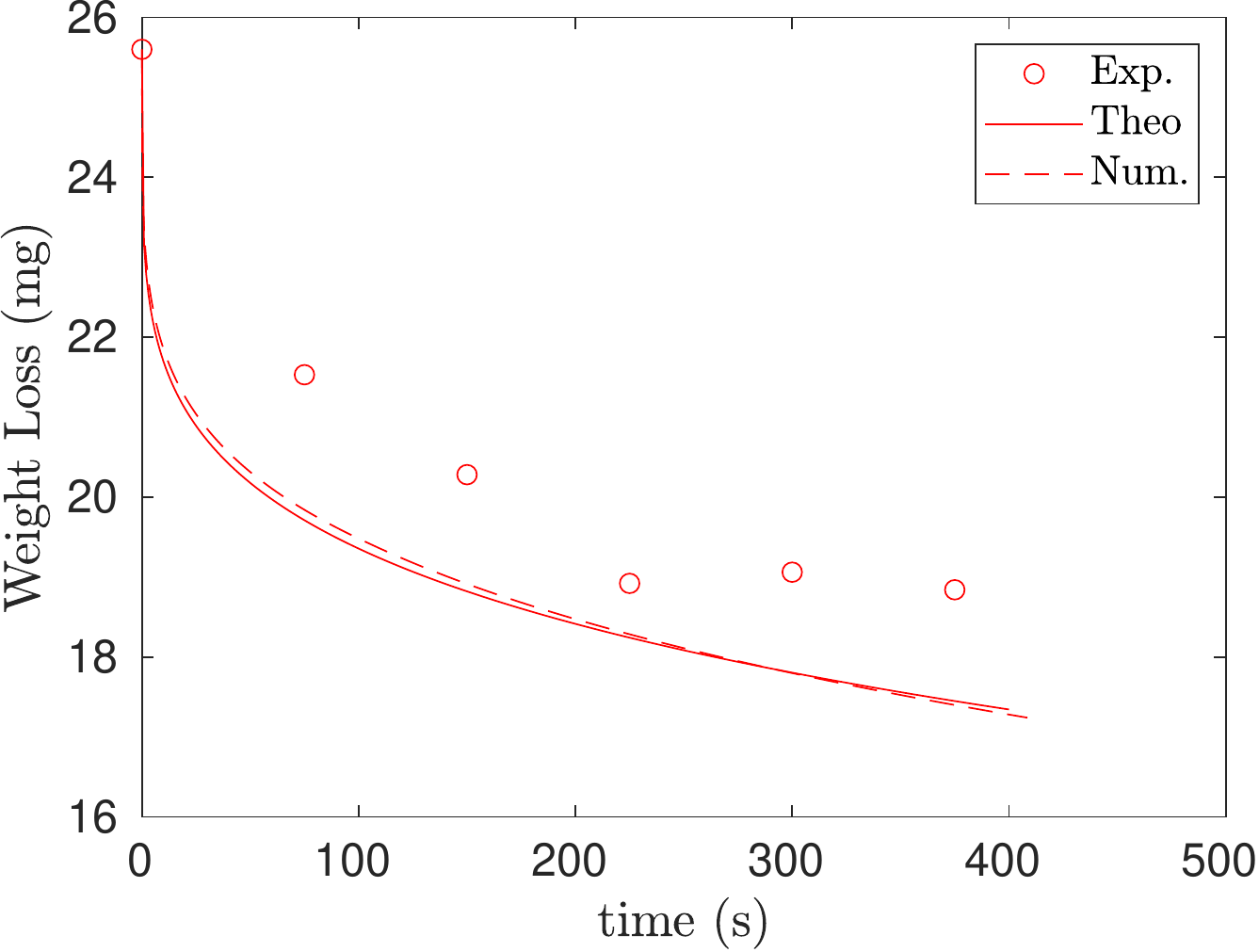}
\caption{TK11 Body Radial}
\label{TK11BR_Wt_AbqExpTh}
\end{subfigure}\\
\begin{subfigure}{.5\textwidth}
\centering
\includegraphics[width=\linewidth]{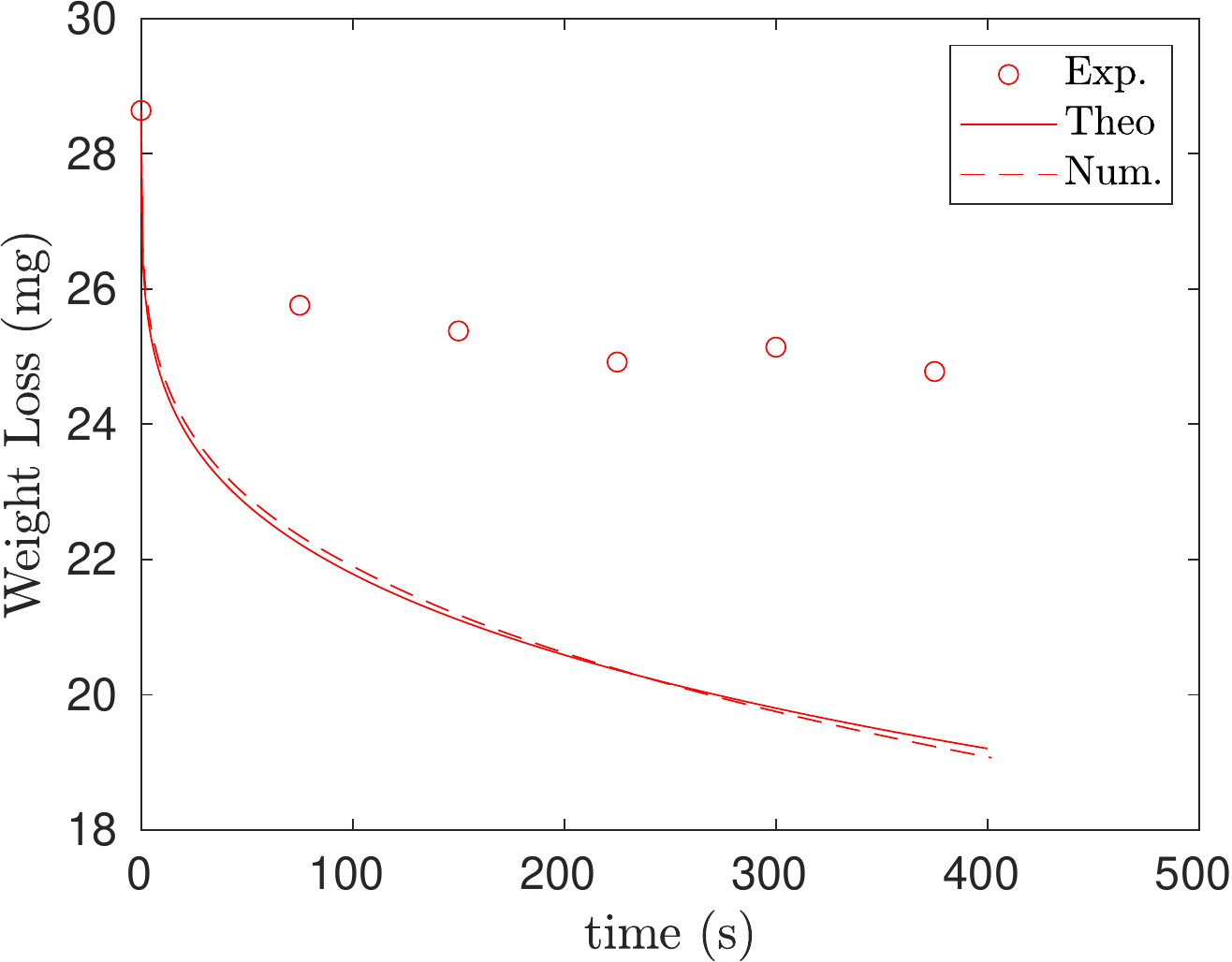}
\caption{TK11 Body Vertical}
\label{TK11BV_Wt_AbqExpTh}
\end{subfigure}%
\begin{subfigure}{.5\textwidth}
\centering
\includegraphics[width=\linewidth]{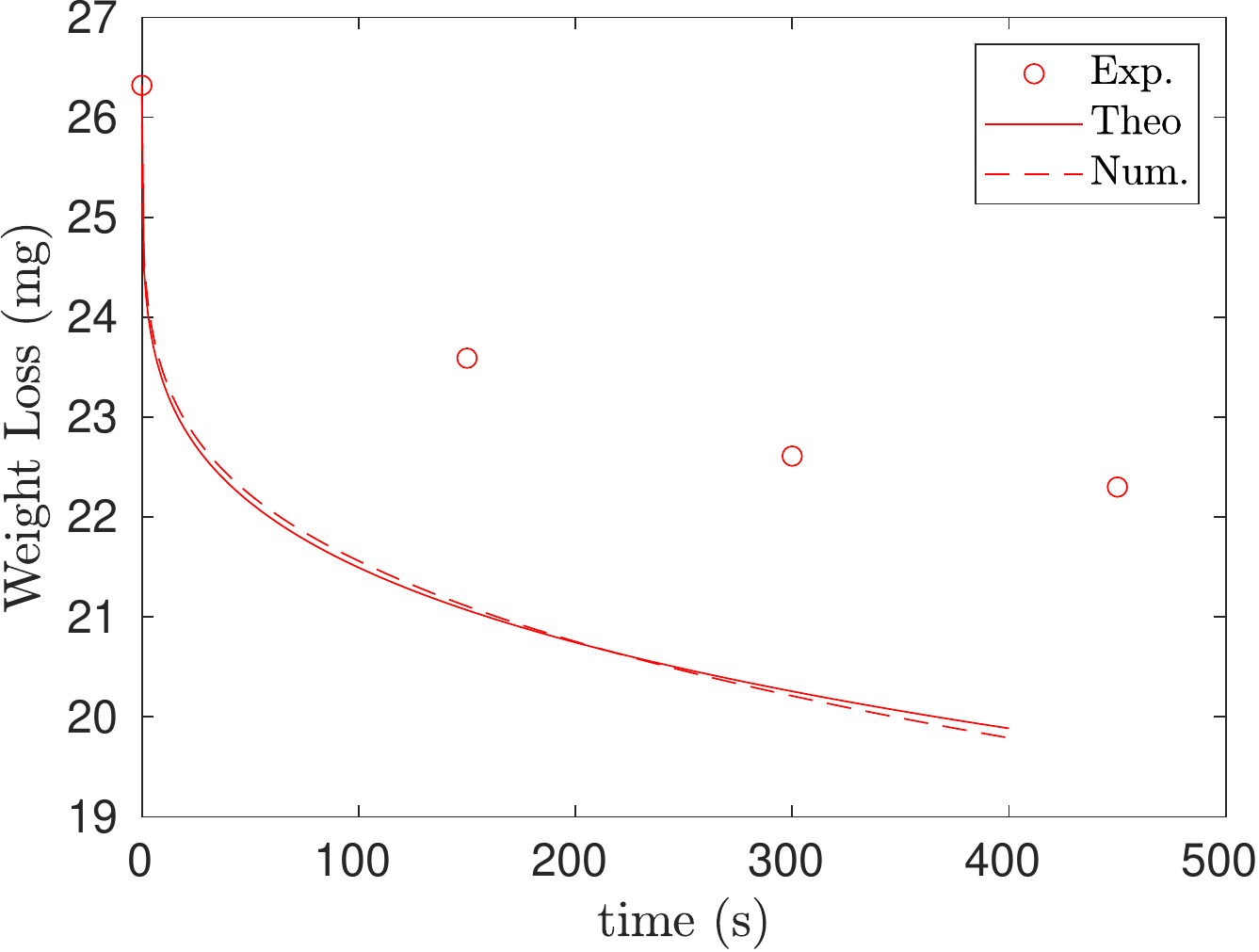}
\caption{TK16 Anterior horn Radial}
\label{TK16AR_Wt_AbqExpTh}
\end{subfigure}%
\caption{Comparison of theoretical weight loss from the sample using parameters obtained from fitting with the experiments (a) Sample taken from the body region of TK11 sample in the circumferential direction, (b) Sample taken from the body region of TK11 sample in the radial direction, (c) Sample taken from the body region of TK11 sample in the vertical direction, (d) Sample taken from the Anthorn region of TK16 sample in the radial direction,}
\label{WtLossCompar2}
\end{figure}

\end{document}